\documentclass[11pt]{amsart}

\usepackage{amssymb} 
\usepackage{hyperref} 
\usepackage{natbib}

\newtheorem{theorem}{Theorem}[section]
\newtheorem{lemma}[theorem]{Lemma}

\newtheorem{proposition}[theorem]{Proposition}
\newtheorem{corollary}[theorem]{Corollary}

\theoremstyle{definition}
\newtheorem{definition}[theorem]{Definition}
\newtheorem{example}[theorem]{Example}

\theoremstyle{remark}
\newtheorem{remark}[theorem]{Remark}

\numberwithin{equation}{section}

\newcommand\ack{\section*{Acknowledgement}}

\begin{document}

\title{On One Number-Theoretic Conception: Towards a New Theory}

\author{\fbox{Rom Varshamov}}
\address{Rom R. Varshamov was a Leading Research Scientist, V.A. Trapeznikov Institute for Control Sciences, Russian Academy of Sciences,
Profsoyuznaya 65, 117997, Moscow, Russia}

\author{Armen Bagdasaryan}
\address{Armen G. Bagdasaryan is a Senior Research Scientist, V.A. Trapeznikov Institute for Control Sciences, Russian Academy of Sciences,
Profsoyuznaya 65, 117997, Moscow, Russia}
\email{abagdasari@hotmail.com, bagdasar@ipu.ru}


\date{September 9, 2009} 



\begin{abstract} \vspace{0.7cm} 
In this paper we present a new mathematical conception based on a new method for ordering the integers. The method relies on the assumption that negative numbers are beyond infinity, which goes back to Wallis and Euler. We also present a new axiom system, the model of which is arithmetics. We define regular method for summation of infinite series which allows us to discover general and unified approach to summation of divergent series, and determine the limits of unbounded and oscillating functions. Several properties for divergent series and explicit formulas for sums of some infinite series are established. A number of finite and new recurrence formulas for Bernoulli numbers are obtained. We rederive some known results, but in a simpler and elementary way, and establish new results by means of techniques of the theoretical background developed. \\ 
\end{abstract}

\maketitle

\begin{center}
\rule{9.5cm}{0.25mm}
\end{center}

\newpage

\tableofcontents 
\addtocontents{toc}{\vspace{20pt}}
\addtocontents{toc}{\baselineskip 20pt}

\baselineskip 18pt

\newpage

\normalsize

\section{Introduction}

\textbf{Historical remarks}. The basic and fundamental concept underlying the foundations of the whole mathematics is the notion of natural number. 
The concept of a negative number appeared much later and emerged in history within the context of algebraic calculations.
Only from the beginning of the sixteenth century onwards, we see the first steps towards negative values, in the form of algebraic terms affected by a negative sign.
It is for the needs of algebra that negative numbers had been introduced to extend natural numbers to the set of all integers. However, 
some properties of negative numbers had been remaining unclear for a long time, 
in particular, the order relation between positive and negative numbers. These gave rise a discussion over a long period of time, in which many of great 
mathematicians had been involved.
Some of them considered negative numbers to be smaller than nothing, while the others believed that negatives were larger than infinity.

So, in seventeenth century there existed at least two fundamentally different approaches to the 
extension of natural numbers to the set of all integers by means of negative numbers: \\
(1) negative numbers are less than nothing (zero), $-1<0$ (Ren\'e Descartes, Albert Girard, Michael Stifel), \\
(2) negative numbers are greater than infinity, $-1>\infty$ (John Wallis, Leonard Euler; and probably Blaise Pascal 
who believed that nothing could be less than nothing \citep{gunt,struik}.

The latter idea, about numbers greater than infinity, comes from treatments of infinitesimals and summation of series. 
John Wallis, in his \textit{Arithmetica Infinitorum} of 1656, advanced the idea that the result of dividing a positive number by a negative one is larger than infinity.
Dealing with the quadrature of curves with equations having negative indices, he came to the statement that the ratio of a positive number to a negative one is greater
than infinity, thus arriving at the conclusion that negative numbers are greater than infinity \citep{Ball, Kline1, Kline2}.

As Wallis, Euler came to the same idea through reasoning about divergent series in his fundamental work \textit{De Seriebus divergentibus} of 1746, presented in 1754 and published in 1760 \citep{eulerseries}. So, Euler also concluded that negative numbers are greater than infinity and $-1>\infty$ \citep{Kline3, sandifer}.

Nevertheless, by the time, when Euler in \textit{``Institutiones Calculi Differentialis''} \citep{euler} proved the rightness of his point of view, 
showing that positive and negative numbers are linked by transition through infinity, the definition of Descartes 
and his followers had already been widely disseminated.
Indeed \citep{gunt}, Isaac Newton, reading lectures at Cambridge from 1673 to 1683, which were edited and published against his wishes in 1707 
by W. Whiston with the title \textit{``Arithmetica Universalis''}, also identified negative numbers as less than nothing, understanding by ``nothing'' zero.
This definition shortly moved to the mathematics textbooks and became widely disseminated in Germany owing to the books of Christian von Wolff. 

\textbf{Subject of research}. This work is motivated by the idea that negative numbers are greater than infinity, expressed by Euler and Wallis. 
Inspired by this hypothesis, that numbers beyond infinity might be negative, we first define a new order relation between negative and positive numbers, 
and develop theoretical statements, which then lead to quite a number of interesting and unexpected results.

The aim of this paper is three-fold:
\begin{itemize}
	\item[(i)] to introduce a new method for ordering the integers where negative numbers are beyond infinity and a new axiom system, 
the model of which is arithmetics, and to study the consequences of this axiomatics;
	\item[(ii)] to develop new mathematical conception and foundations of a new theory on the basis of (i);
	\item[(iii)] using (ii) and applying it to concrete topics in different areas of mathematics, to show the possibilities of a new theoretical background
				by reproving and/or rederiving known results, using a simpler approach, and by obtaining new results by means of new techniques; these results can not be explained by and deduced from the existing theories and at the same they do not contradict them.
\end{itemize}

The starting points of the work are:
\begin{itemize}
	\item[(i)] a new method for ordering the integers, from which we get  $\mathbb{Z}=[0,1,2,...,-2,-1]$;
	\item[(ii)] a new class of real regular functions $f(\cdot)$ defined on $\mathbb{Z}$ and the definition of $\sum_{a}^{b}f(\cdot)$ that extends 
the classical definition to the case $b<a$;
	\item[(iii)] a set of axioms imposed on regular functions; these axioms define a new regular method for summation of infinite series. 
\end{itemize}

Thus, we find a unified approach to summation of divergent series, and  
to determination of limits of unbounded and oscillating functions.

The new theoretical background developed, being consistent with what has been obtained 
in the framework and on the basis of existing theories with the help of the known means of analysis, considerably 
extends the possibilities of modern mathematics from general standpoints, and leads to the discoveries of a qualitatively new nature.

The paper is organized as follows. In Section \ref{sec:intro} we give basic definitions, introduce new axiom system, and obtain some propositions resulting from the system; 
in Section \ref{sec:main} we give a number of propositions concerning infinite series, including several properties of divergent series,
and some propositions on the limits of unbounded and oscillating functions;
Section \ref{sec:bern} deals, mainly, with the Bernoulli numbers and with some results around them; 
Section \ref{sec:trig} concerns the limits of trigonometric functions and the sums of trigonometric series; 
in Section \ref{sec:exten} we show how the class of regular functions can be substantially extended if constructed with use of non-elementary functions.
We conclude in Section \ref{sec:conc} with discussion of the obtained results and with some aspects of future research
which will be based on our theory formulated as a paradigm.

\textbf{Notation}. Throughout the paper, we use the following notation. We use the symbols
$$
\mathbb{Z}, \mathbb{N}, \mathbb{R}
$$
to denote the set of integers, the set of the natural numbers, and the set of real numbers, respectively. 
And we denote by $a,b,c\in\mathbb{Z}$ the integer numbers, and by $n,m,k\in\mathbb{N}$ the natural numbers.

\vspace{1cm}
\begin{center}
\rule{9.5cm}{0.25mm}
\end{center}

\newpage

\section{New Ordering of Integers, Axiom System, and Regular Functions}\label{sec:intro}

In this section we give basic definitions and concepts of new axiomatics,
the model of which is arithmetics, and also present main propositions resulting from the axiom system. 

\subsection{New Method for Ordering the Integers and Axiom System}

Consider the set of all integer numbers $\mathbb{Z}$. We introduce a new method for ordering the integers (first introduced, but in other form, in \citep{D}, 
and then in \citep{bag}) as follows.
\begin{definition} 
We shall say that $a$ precedes $b$, $a, b \in \mathbb{Z}$, and write $a\prec b$, if the inequality $\frac{-1}{a}<\frac{-1}{b}$ holds; $a\prec b \Leftrightarrow \frac{-1}{a}<\frac{-1}{b}$ \footnote{assuming by convention $0^{-1}=\infty$}.
\end{definition}
From this method of ordering it follows that any positive integer number, including zero, precedes any negative integer number, and the set $\mathbb{Z}$ has zero as the first element and $-1$ as the last element, i. e. we have $\mathbb{Z}=[0, 1, 2,...-2, -1]$\footnote{the set $\mathbb{Z}$ can be
geometrically represented as cyclically closed}. And the following two essential axioms of order hold:
\begin{itemize}
\item if $a\prec b$ and $b\prec c$ then $a\prec c$ (transitivity)
\item if $a\neq b$ then either $a\prec b$ or $b\prec a$ (connectedness)
\end{itemize}
\begin{definition}
The function $f(x)$, $x\in \mathbb{Z}$ , is called \textit{regular} if there exists an elementary \footnote{the function which is determined by formulas containing a finite number of algebraic or trigonometric operations performed over argument, function and constants \citep{N}} function $F(x)$ such that $F(z+1)-F(z)=f(z), \; \forall z\in \mathbb{Z}$. The function $F(x)$ is said to be a \textit{generating function} for $f(x)$.
\end{definition}
\begin{remark}
If $F(x)$ is a generating function for $f(x)$, then the function $F(x)+C$, where $C$ is a constant, is also a generating function for $f(x)$. So, any function $F(x)$ which is generating for $f(x)$ can be represented in the form $F(x)+C(x)$, where $C(x)$ is a periodic function with the period $1$. It is clear that \textit{any linear combination of regular functions is also a regular function}.
\end{remark}

Suppose $f(x)$ is a function of real variable defined on $\mathbb{Z}$ and $\mathbb{Z}_{a, b}$ is a part of $\mathbb{Z}$ such that 
$$ \mathbb{Z}_{a, b} = \left\{              
									\begin{array}{ll}                   
									[a, b] & (a\preceq b)\\                   
									\mathbb{Z}\setminus (b, a) & (a\succ b)              
									\end{array}       
						 \right.
$$
where $\mathbb{Z}\setminus (b, a)=[a, -1]\cup[0, b]$ \footnote{the segment $[a, b]$ (interval $(a, b)$), as usual, is the set of all integers $x$ such that $a\preceq x\preceq b$ ($a\prec x\prec b$)}.

\begin{definition}
For any $a, b\in \mathbb{Z}$
\begin{equation} 
\sum_{u=a}^b{f(u)}=\sum_{u\in \mathbb{Z}_{a, b}}{f(u)}\footnote{observing the established order of elements $\mathbb{Z}_{a, b}$}. \label{def:first_sum_eq}
\end{equation}
\end{definition}
This definition satisfies the condition of generality and has a real sense for any integer values of $a$ and $b$ ($a ^{>}_{<} b$).
As it is known a sum $\sum_a^b$ is not defined for $b<a$, it is generally understood to be zero.
Our definition of $\sum_{a}^{b}f(n)$ extends  the classical definition to the case $b<a$. 
With this definition  a sum $\sum_a^bf(u)$ is defined for arbitrary limits of summation ($a ^{>}_{<} b$) and in such a way 
that its functional dependence for $a\leq b$ retains its analytical expression (i. e. the method is regular), $F(a, b)=\sum_{u=a}^b {f(u), \; (a ^> _< b)}$.

\begin{remark}
The equality analogous to (\ref{def:first_sum_eq}) can also be defined for products $\prod_a^b$ with arbitrary limits $a$ and $b$, $(a ^> _< b)$.
\end{remark}

The set $\mathbb{Z}_{a, b}$, depending on the elements $a$ and $b$, can be both finite and infinite. Thus, the sum on the right-hand side of (2.1) may become an infinite series. 
For this reason, we have to assign a numerical meaning to $\sum$
\footnote{the notions of a sum and a limit (see definition \ref{def:limit_seq}) are understood in a usual, classical, sense, $\sum_{n=\alpha}^{\alpha+k}f(n)=f(\alpha)+f(\alpha+1)+f(\alpha+2)+...+f(\alpha+k), \alpha\in \mathbb{Z}, k\in \mathbb{N}$, so we do not introduce new symbols for $\sum$ and $\lim$}; Following L. Euler\footnote{who was convinced that ``to every series could be assigned a number'', which is the first Euler's principle on infinite series}, 
we take as a postulate the assertion that \textit{any series $\sum_{u=1}^{\infty}f(u)$,
where $f(u)$ is a regular function, has a certain numeric value}.

We also introduce and require that regular functions satisfy the following, quite natural, system of axioms.
\begin{enumerate}
\item If $S_n=\sum_{u=a}^n{f(u)} \;\; \forall n$, then $\lim_{n\rightarrow \infty}S_n=\sum_{u=a}^\infty {f(u)}$\footnote{$n\rightarrow \infty$ means that $n$ unboundedly increases, without changing the sign}. \label{axiom:1}    \\
\item If $S_n=\sum_{u=1}^{n/2}{f(u)} \;\; \forall n$, then $\lim_{n\rightarrow \infty}S_n=\sum_{u=1}^\infty {f(u)}$. \label{axiom:2}   \\
\item If $\sum_{u=a}^\infty {f(u)}=S$, then $\sum_{u=a}^\infty {\alpha f(u)}=\alpha S, \; \alpha\in \mathbb{R}$. \label{axiom:3}     \\
\item  If $\sum_{u=a}^\infty {f_1(u)}=S_1$ and $\sum_{u=a}^\infty {f_2(u)}=S_2$, then  $\sum_{u=a}^\infty {(f_1(u)+f_2(u))}= S_1+S_2$.  \label{axiom:4}  \\
\item For any $a$ and $b$, $a\leq b$: $F(b+1)-F(a)=\sum_{u=a}^b{f(u)}$. \label{axiom:5}   \\
\item If $G=[a_1, b_1]\cup[a_2, b_2]$, $[a_1, b_1]\cap[a_2, b_2]=\emptyset$, then $\sum_{u\in G}{f(u)}=\sum_{u=a_1}^{b_1}{f(u)}+\sum_{u=a_2}^{b_2}{f(u)}$. \label{axiom:6} 
\end{enumerate}

The axioms (\ref{axiom:1})--(\ref{axiom:6}) define the method of summation of infinite series, which is regular, since, due to axiom (\ref{axiom:1}), it sums every convergent series to its usual sum in a classical sense. 
\begin{definition} \label{def:limit_seq}
A number $A$ is said to be the \textit{limit} of a numeric sequence $F(1)$, $F(2)$,..., $F(n)$,... (function of integer argument), i. e. $\lim_{n\rightarrow\infty}F(n)=A$, if $\sum_{u=1}^{\infty}f(u)=A$, where $F(1)=f(1)$ and $F(u)-F(u-1)=f(u) \; \; (u>1)$, i. e.
$$
F(1)+(F(2)-F(1))+...+(F(n)-F(n-1))+...=A. 
$$
\end{definition}

\begin{remark}
Note that the Definition \ref{def:limit_seq} in particular coincides with the classical definition of limit of sequences, which is convergent in a usual sense, and reduces the problem of existence of the limit of functions of integer argument $F(n)$ to finding the sum of the series
$$
F(1)+(F(2)-F(1))+...+(F(n)-F(n-1))+...
$$
This, in turn, allows us to claim (in force of our postulate) that any elementary function of integer argument defined on $\mathbb{Z}$ has a certain limit. 
\end{remark}

Relying on axiom (\ref{axiom:1}), the following two rules of limit are easily derived:
\begin{enumerate}
	\item The limit of constant equals to the same constant, i. e. $$\lim_{n\rightarrow\infty}A=A$$  \label{rule1}
	\item The limit of algebraic sum of finite number of sequences equals to the algebraic sum of limits of the sequences, i. e. $$\lim_{n\rightarrow\infty}\sum_{u=1}^k{\alpha_u F_u(n)}=\sum_{u=1}^k{\alpha_u\lim_{n\rightarrow\infty}F_u(n)}$$  \label{rule2}
	where $\alpha_u$ are real numbers. \\

Other two rules - limit of product and limit of quotient - are, generally speaking, true only for convergent sequences. However, further in this work, we encounter with limits of  sums of infinite terms. So, without going into details, we postulate the rule \\
	\item The limit of the sum of convergent series equals to the sum of limits of series' terms, \label{rule3}
	e. g. $$\lim_{n\rightarrow\infty}\sin n\theta=\theta\lim_{n\rightarrow\infty}n-\frac{\theta^3}{3!}\lim_{n\rightarrow\infty}n^3+\frac{\theta^5}{5!}\lim_{n\rightarrow\infty}n^5-... \; \; \; \; \; |x|<\infty$$  
	$$\lim_{n\rightarrow\infty}\cos n\theta=1-\frac{\theta^2}{2!}\lim_{n\rightarrow\infty}n^2+\frac{\theta^4}{4!}\lim_{n\rightarrow\infty}n^4-... \; \; \; \; \; |x|<\infty$$
	$$\lim_{n\rightarrow\infty}\sum_{u=1}^n{(-1)^{u-1}\frac{\sin ux}{u^{2k+1}}}=\sum_{u=1}^\infty{(-1)^{u-1}\frac{\sin ux}{u^{2k+1}}}=x\lim_{n\rightarrow\infty}p_0(n)-\frac{x^3}{3!}\lim_{n\rightarrow\infty}p_1(n)+... \; \; \forall n$$
	$$\lim_{n\rightarrow\infty}\sum_{u=1}^n{(-1)^{u-1}\frac{\cos ux}{u^{2k}}}=\sum_{u=1}^\infty{(-1)^{u-1}\frac{\cos ux}{u^{2k}}}=\lim_{n\rightarrow\infty}p_0(n)-\frac{x^2}{2!}\lim_{n\rightarrow\infty}p_1(n)+... \; \; \forall n$$ 
	where $p_i(n)=\sum_{u=1}^n\frac{(-1)^{u-1}}{u^{2(k-i)}}$. \\
\end{enumerate}

\begin{remark} \label{rem:a_bigger_b}
If $a\succ b$ then, by definition, $\mathbb{Z}_{a, b}=[a, -1]\cup[0, b]$, $[a, -1]\cap[0, b]=\emptyset$ and, in view of axiom (\ref{axiom:6}), for any regular function $f(x)$
\begin{equation}
\sum_{u=a}^b f(u)=\sum_{u=a}^{-1}f(u)+\sum_{u=0}^b f(u), \; \; \; a\succ b. \label{eq:a_bigger_b}
\end{equation}
\end{remark}

\begin{remark} \label{rem:sum_of_two_intervals}
Using axiom (\ref{axiom:6}),  it can be easily shown that if $G=\mathbb{Z}_{a_1, b_1}\cup \mathbb{Z}_{a_2, b_2}$ and $\mathbb{Z}_{a_1, b_1}\cap \mathbb{Z}_{a_2, b_2}=\emptyset$ then $$\sum_{u\in G}f(u)=\sum_{u=a_1}^{b_1}f(u)+\sum_{u=a_2}^{b_2}f(u)$$
\end{remark}

\subsection{Propositions Involving Regular Functions and Resulting From New Axiom System}

\begin{proposition} \label{prop:a_a-1}
If $f(x)$ is a regular function and $a\in \mathbb{Z}$ is a fixed number, then 
\begin{equation}
\sum_{u=a}^{a-1}f(u)=\sum_{u\in \mathbb{Z}}f(u)
\end{equation}
\end{proposition}
\begin{proof}
Let us consider two cases. \\
Case 1. If $a=0$ then $0\prec -1$ and, by definition, $\mathbb{Z}_{0,-1}=[0, 1, 2,..., -2, -1]=\mathbb{Z}$. Hence, 
$$\sum_{u=a}^{a-1}f(u)=\sum_{u\in \mathbb{Z}}f(u) \;\;\;\;\; (a=0)$$
Case 2. If $a\neq 0$ then $a\succ a-1$ and, by definition, $\mathbb{Z}_{a, a-1}=[a, -1]\cup [0, a-1]=[0, a-1]\cup [a, -1]=[0, 1, 2,..., -2, -1]$, 
i. e. $\mathbb{Z}_{a, a-1}=\mathbb{Z}$, and we get  
$$\sum_{u=a}^{a-1}f(u)=\sum_{u\in \mathbb{Z}}f(u) \;\;\;\;\; (a\neq 0)$$
This concludes the proof.
\end{proof}

\begin{proposition} \label{prop:pair_of_numbers}
For any numbers $m$ and $n$ such that $m\prec n$
\begin{equation}
\sum_{u=m}^{n}f(u)=\sum_{u=-n}^{-m}f(-u) \label{eq:pair_of_numbers}
\end{equation}
\end{proposition}
\begin{proof}
The proof is based on the fact that the order relation between pairs of numbers $m, n$ and $-n, -m$ is the same.
\end{proof}

\begin{proposition} \label{prop:mean}
Let $f(x)$ be a regular function and let $a$, $b$, $c$ be any integer numbers such that $b\in \mathbb{Z}_{a, c}$. Then 
\begin{equation}
\sum_{u=a}^{c}f(u)=\sum_{u=a}^{b}f(u)+\sum_{u=b+1}^{c}\!\!\acute{}\;f(u),  \label{eq:prop_mean}
\end{equation}
where the prime on the summation sign means that $\sum\limits_{u=b+1}^{c}\!\!\!\!\acute{}\;f(u)=0$ for $b=c$.
\end{proposition}
\begin{proof}
For $c=b$ or $c=a$, there is nothing to prove. So, let $c\neq b$ and $c\neq a$. The proof consists of two steps.
\begin{itemize}
\item[(1)] If $a\prec c$ then $\mathbb{Z}_{a, c}=[a, c]$, i. e. $[a, c]=[a, b]\cup [b+1, c]$, where $[a, c]=\\=[a, b]\cap [b+1, c]=\emptyset$. In view of axiom (\ref{axiom:6}), we obtain
$$\sum_{u=a}^{c}f(u)=\sum_{u=a}^{b}f(u)+\sum_{u=b+1}^{c}f(u) \;\;\;\;\; (b\neq c)$$

\item[(2)] If $a\succ c$ then $\mathbb{Z}_{a, c}=[a, -1]\cup [0, c]$ and, therefore, there are two cases: either $b\in [a, -1]$ or $b\in [0, c]$.
\begin{itemize}
\item[(a)] If $b\in [a, -1]$ then $a\preceq b\preceq -1$, and $[a, -1]=[a, b]\cup [b+1, -1]$, assuming that $[b+1, -1]=\emptyset$ for $b=-1$. Hence we have 
$\mathbb{Z}_{a, c}=[a, -1]\cup [0, c]=\\=[a, b]\cup\big([b+1, -1]\cup[0, c]\big)=\mathbb{Z}_{a, b}\cup \mathbb{Z}_{b+1, c}$, 
where $\mathbb{Z}_{a, b}\cap \mathbb{Z}_{b+1, c}=\emptyset$. 
According to the Remark \ref{rem:sum_of_two_intervals}, we get
$$\sum_{u=a}^{c}f(u)=\sum_{u=a}^{b}f(u)+\sum_{u=b+1}^{c}f(u) \;\;\;\;\; (b\neq c)$$
\item[(b)] If $b\in [0, c]$ then $0\preceq b\prec c$, and $[0, c]=[0, b]\cup [b+1, c]$, i. e. 
$\mathbb{Z}_{a, c}=\\=[a, -1]\cup[0, b]\cup[b+1, c]=\mathbb{Z}_{a, b}\cup \mathbb{Z}_{b+1, c}$, where $\mathbb{Z}_{a, b}\cap \mathbb{Z}_{b+1, c}=\emptyset$. Hence we again get
$$\sum_{u=a}^{c}f(u)=\sum_{u=a}^{b}f(u)+\sum_{u=b+1}^{c}f(u) \;\;\;\;\; (b\neq c)$$
\end{itemize}
\end{itemize}
This concludes the proof.
\end{proof}

\begin{corollary}[Proposition \ref{prop:mean}] 
For any regular function $f(x)$
\begin{equation}
\sum_{u\in \mathbb{Z}_{-n,n}}f(u)=\sum_{u=0}^{n}f(u)+\sum_{u=1}^{n}f(-u) \;\;\;\;\; \forall{n}  \label{eq:cor_prop_mean}
\end{equation}
\end{corollary}
\begin{proof}
This equality immediately follows from the formula (\ref{def:first_sum_eq}) (letting $a=-n, b=n$), axiom (\ref{axiom:6}), and Proposition \ref{prop:pair_of_numbers}. Indeed, substituting $-n$ and $n$ for $a$ and $b$ correspondingly in formula (\ref{def:first_sum_eq}), and having that $\mathbb{Z}_{-n, n}=[-n, -1]\cup [0, n]$ and $[-n, -1]\cap [0, n]=\emptyset$, we, in view of axiom (\ref{axiom:6}), get $$\sum_{u\in \mathbb{Z}_{-n,n}}f(u)=\sum_{u=0}^{n}f(u)+\sum_{u=-n}^{-1}f(u)$$
Now using Proposition \ref{prop:pair_of_numbers}, we finally obtain
$$\sum_{u\in \mathbb{Z}_{-n,n}}f(u)=\sum_{u=0}^{n}f(u)+\sum_{u=1}^{n}f(-u)$$
\end{proof}

\begin{proposition}
Suppose $f(x)$ is a regular function. Then
\begin{equation}
\sum_{u=a}^{a-1}f(u)=0 \;\;\;\;\; \forall{a\in \mathbb{Z}}  \label{eq:sum_a_a-1_zero}
\end{equation}
or, which is the same in view of Proposition \ref{prop:a_a-1},
\begin{equation}
\sum_{u\in \mathbb{Z}}f(u)=0 \label{eq:sum_zero}
\end{equation}
\end{proposition}
\begin{proof}
The proof is in two steps.
\begin{itemize}
\item[(1)] Let $a=0$. Then $a\prec a-1$, and according to axiom (\ref{axiom:5})
$$\sum_{u=a}^{a-1}f(u)=F(a)-F(a)=0$$
i. e.
$$\sum_{u=a}^{a-1}f(u)=0 \;\;\;\;\; (a=0)$$
\item[(2)] If $a\neq 0$ then $a\succ a-1$ and taking into account Remark \ref{rem:a_bigger_b} and axiom \ref{axiom:5}, we obtain 
$$\sum_{u=a}^{a-1}f(u)=\sum_{u=a}^{-1}f(u)+\sum_{u=0}^{a-1}f(u)=F(0)-F(a)+F(a)-F(0)=0$$
and
$$\sum_{u=a}^{a-1}f(u)=0 \;\;\;\;\; (a\neq 0)$$
\end{itemize}
The proof is completed.
\end{proof}

The following two lemmas and their corollaries demonstrate some properties of new number axis and shed a light on its structure.

\begin{lemma} \label{lemma:1} 
The following equality holds
\begin{equation}
\lim_{n\rightarrow\infty}(n, -n)=\emptyset  \label{eq:1_lemma1}
\end{equation}
\end{lemma}
\begin{proof}
Since $-n\succ n$, then, by definition, $\mathbb{Z}_{-n, n}=[-n, -1]\cup [0, n]$. Therefore, for arbitrary number $a$ there exists $n_0=|a|$ such that for every natural number $n\geq n_0$ the relation $a\in \mathbb{Z}_{-n, n}$ holds. Whence it follows that
\begin{equation}
\lim_{n\rightarrow\infty}\mathbb{Z}_{-n, n}=\mathbb{Z} \label{eq:lemma1}
\end{equation}
But $\mathbb{Z}_{-n, n}=\mathbb{Z}\setminus (n, -n)$. Therefore, passing to the limit and taking into account (\ref{eq:lemma1}), we obtain $$\lim_{n\rightarrow\infty}(n, -n)=\emptyset$$
This completes the proof.
\end{proof}
Relying on (\ref{eq:lemma1}) and using the analogy of correspondence of partial sums to infinite series, it can be easily shown that
\begin{equation}
\lim_{n\rightarrow\infty}\sum_{u\in \mathbb{Z}_{-n, n}}f(u)=\sum_{u\in \mathbb{Z}}f(u) \label{eq:partial_sum}
\end{equation}

As we already know from formula (\ref{eq:sum_zero}) $\sum_{u\in \mathbb{Z}}f(u)=0$. Then using (\ref{eq:cor_prop_mean}), (\ref{eq:partial_sum}) and axiom (\ref{axiom:1}), we have
\begin{equation}
\sum_{u=0}^{\infty}f(u)+\sum_{u=1}^{\infty}f(-u)=0  \label{eq:1_from_lemma1}
\end{equation}

\begin{lemma} \label{lemma2} 
For any elementary function $F(x)$ defined on $\mathbb{Z}$
\begin{equation}
\lim_{n\rightarrow\infty}F(n+1+a)=\lim_{n\rightarrow\infty}F(-n+a) \;\;\;\;\; \forall a   \label{eq:lemma2}
\end{equation}
\end{lemma}
\begin{proof}
Let $F(x+1+a)-F(x+a)=f(x)$. According to (\ref{eq:a_bigger_b}), for any natural $n$
$$\sum_{u\in \mathbb{Z}_{-n, n}}f(u)=\sum_{u=-n}^{-1}f(u)+\sum_{u=0}^{n}f(u)$$
Hence, in view of axiom (\ref{axiom:5})
$$\sum_{u\in \mathbb{Z}_{-n, n}}f(u)=F(a)-F(-n+a)+F(n+1+a)-F(a)=F(n+1+a)-F(-n+a)$$
By definition, the function $f(x)$ is regular. Thus, passing to the limit, $n\rightarrow\infty$, and taking into account (\ref{eq:sum_zero}) and (\ref{eq:partial_sum}), we obtain
$$\lim_{n\rightarrow\infty}F(n+1+a)=\lim_{n\rightarrow\infty}F(-n+a),$$
which completes the proof.
\end{proof}

\begin{corollary}[Lemma \ref{lemma2}]
For any elementary function $F(x)$ defined on $\mathbb{Z}$
$$\lim_{n\rightarrow\infty}F(-(n+1))=\lim_{n\rightarrow\infty}F(n).$$
\end{corollary}
In particular,
\begin{equation}
\lim_{n\rightarrow\infty}(-(n+1))=\lim_{n\rightarrow\infty}n  \label{eq:1_cor_lemma2}
\end{equation}
\begin{equation}
\lim_{n\rightarrow\infty}(n+1)=\lim_{n\rightarrow\infty}(-n)   \label{eq:2_cor_lemma2}
\end{equation}
i. e.
$$\lim_{n\rightarrow\infty}n=-\frac{1}{2}.$$
Since $n=\sum_{u=1}^{n}1 \;\;\; \forall n$, then, using axiom (\ref{axiom:1}), $\lim_{n\rightarrow\infty}n=\sum_{u=1}^{\infty}1$
and
\begin{equation}
\sum_{u=1}^{\infty}1=-\frac{1}{2}.  \label{eq:sum_of_1}
\end{equation}

\begin{corollary}[Lemma \ref{lemma2}]
For any odd function $\psi(x)$ defined on $\mathbb{Z}$, according to formula (\ref{eq:lemma2}) for $a=-1/2$
\begin{equation}
\lim_{n\rightarrow\infty}\psi(n+\frac{1}{2})=0.  \label{eq:3_cor_lemma2}
\end{equation}
\end{corollary}

\begin{lemma} \label{lemma3}      
For any regular function $f(x)$
\begin{equation}
\sum_{u=a}^{b}f(u)=-\sum_{u=b+1}^{a-1}f(u) \;\;\; \forall a, b\in \mathbb{Z}    \label{eq:lemma3}
\end{equation}
\end{lemma}
\begin{proof}
Letting $c=a-1$ in (\ref{eq:prop_mean}) and taking into account (\ref{eq:sum_a_a-1_zero}), we immediately get (\ref{eq:lemma3}).
\end{proof}

From formula (\ref{eq:lemma3}), using (\ref{eq:pair_of_numbers}) and (\ref{eq:1_cor_lemma2}), we obtain

\begin{equation}
\sum_{u=1}^{-n}f(u)=-\sum_{u=0}^{n-1}f(-u)  \label{eq:1_from_lemma3}
\end{equation}
and
\begin{equation}
\sum_{u=1}^{\infty}f(u)=-\sum_{u=0}^{\infty}f(-u)   \label{eq:2_from_lemma3}
\end{equation}
Indeed, putting $a=0$ and $b=-n$ in (\ref{eq:lemma3}), we have
$$
\sum_{u=0}^{-n}f(u)=-\sum_{u=-n+1}^{-1}f(u)
$$
and using (\ref{eq:pair_of_numbers}), we get
$$
\sum_{u=0}^{-n}f(u)=-\sum_{u=1}^{n-1}f(-u)
$$
and, finally,
$$
\sum_{u=1}^{-n}f(u)=-\sum_{u=0}^{n-1}f(-u)
$$
i.e. the equality (\ref{eq:1_from_lemma3}).

Now, substituting $n$ by $n+1$ in (\ref{eq:1_from_lemma3}) and passing to the limit with regard to (\ref{eq:1_cor_lemma2}), we obtain (\ref{eq:2_from_lemma3}). The formula (\ref{eq:2_from_lemma3}) is, in fact, the equality (\ref{eq:1_from_lemma1}).

The same line of reasoning allows us to find analogous formulas for products. For example, we have
\begin{proposition} \label{prop:product}
Suppose $f(x)$ is a regular function defined on $\mathbb{Z}$ such that $\prod_{u=1}^{0}f(u)=1$. Then
\begin{equation}
\prod_{u=0}^{-n}f(u)=\left(\prod_{u=1}^{n-1}f(-u)\right)^{-1}  \label{eq:product}
\end{equation}
\end{proposition}

We give some consequences resulting from (\ref{eq:1_from_lemma3}) and (\ref{eq:product}). From (\ref{eq:1_from_lemma3}) it follows that all odd Bernoulli numbers (except for $B_1=1/2$) are equal to zero. Let us take
$$B_k(n)=\frac{1}{k+1}\sum_{u=0}^k \biggl(\begin{matrix}k+1\\ u\end{matrix}\biggr)B_u n^{k+1-u}=\sum_{u=1}^{n}u^k \;\; -$$
the Bernoulli polynomial.

Since $B_k(n)-B_k(n-1)=n^k$, the function $f(x)=x^k$ is regular. Then, in view of (\ref{eq:1_from_lemma3}), we have
\begin{equation}
B_k(-n)=(-1)^{k-1}B_k(n-1) \label{eq:bernoulli}
\end{equation}
Hence,
$$
B_k(n)-(-1)^{k-1}B_k(-n)=\frac{2}{k+1}\sum_{u=0}^{\left[(k-1)/2\right]}\biggl(\begin{matrix}k+1\\ 2u+1\end{matrix}\biggr)B_{2u+1}n^{k-2u}=n^k,
$$
from which we get that $B_1=1/2$ and $B_{2u+1}=0, \; u=1,2,3,...$.

Now, relying on (\ref{eq:product}), we extend the function $n!$, defined for natural $n$, to the set of all integers. Let $\lambda(s)=\sum_{u=1}^{s}u=s! \; \forall s$. Then, formally, $\lambda(0)=\sum_{u=1}^{0}u=0!=1$. This means that the function $f(x)=x$ satisfies the Proposition \ref{prop:product}, which allows us to apply the formula (\ref{eq:product}) to this function. Putting $f(u)=u$ in (\ref{eq:product}), one can see that the function $\lambda(s)$, in analogy with the gamma-function, has simple poles at $s=-n$ with residue $(-1)^{n-1}\frac{1}{(n-1)!}$. Taking this into account, we in particular get
$$
\biggl(\begin{matrix}n\\ -m\end{matrix}\biggr)=0
$$
\begin{equation}
\biggl(\begin{matrix}-n\\ m\end{matrix}\biggr)=(-1)^{m}\biggl(\begin{matrix}n+m-1\\ m\end{matrix}\biggr) \;\;\;\;\;\;\;\;\;\;\;\;\;\;\;\;\; \textup{\citep{M}}  \label{eq:combi}
\end{equation}
$$
\biggl(\begin{matrix}-n\\ -m\end{matrix}\biggr)=(-1)^{n+m}\biggl(\begin{matrix}m-1\\ n-1\end{matrix}\biggr) \;\;\;\;\;\;\;\;\;\;\;\;\;\;\;\;\;\; \textup{\citep{M}}
$$

\vspace{1cm}
\begin{center}
\rule{9.5cm}{0.25mm}
\end{center}

\newpage

\section{Functions, Sequences, and Series}\label{sec:main}

In this section we give a number of theorems regarding limits of unbounded and oscillating functions, 
infinite series and their properties.  

\subsection{Series and Limits of Functions}

Using (\ref{eq:2_from_lemma3}) we obtain
\begin{theorem} \label{theorem:inf_sum}
For any even regular function $f(x)$
\begin{equation}
\sum_{u=1}^{\infty}f(u)=-\frac{f(0)}{2}  \label{eq:inf_sum}
\end{equation}
independently on whether the series is convergent or not in a usual sense.
\end{theorem}
\begin{proof}
Theorem will be proved by two ways.
\begin{itemize}
\item[(1)] Using (\ref{eq:2_from_lemma3}), we get 
\begin{eqnarray}
&&\sum\limits_{u=1}^{\infty}f(u) = -\sum\limits_{u=0}^{\infty}f(-u) \Rightarrow\ \sum\limits_{u=1}^{\infty}f(u)=-\sum\limits_{u=0}^{\infty}f(u)\Rightarrow \sum\limits_{u=1}^{\infty}f(u) = -f(0)-\sum\limits_{u=1}^{\infty}f(u) \Rightarrow \nonumber \\ 
&&\Rightarrow \nonumber
2\sum\limits_{u=1}^{\infty}f(u)=-f(0) \Rightarrow \sum\limits_{u=1}^{\infty}f(u)=-f(0)/2.
\end{eqnarray}
\item[(2)] Using (\ref{eq:1_from_lemma1}), we get
\begin{eqnarray}
&&\sum\limits_{u=0}^{\infty}f(u)+\sum\limits_{u=1}^{\infty}f(-u) =0 \Rightarrow f(0)+\sum\limits_{u=1}^{\infty}f(u)+ \sum\limits_{u=1}^{\infty}f(u)=0\Rightarrow
2\sum\limits_{u=1}^{\infty}f(u) = \nonumber \\
&&= \nonumber
-f(0)\Rightarrow \sum\limits_{u=1}^{\infty}f(u)=-f(0)/2.
\end{eqnarray}
\end{itemize}
\end{proof}

\begin{example}
We consider some examples of both convergent and divergent series.
\begin{itemize}
\item[(1)] Convergent series 
$$
\sum_{u=1}^{\infty}\frac{1}{4u^2-1}=\frac{1}{2} \;\;\;\;\;\;\;\;\;\;\;\;\;\;\;\;\;\;\;\;\;\;\;\;\;\;\;\;\;\;\;\;\; \textup{\citep{I}}
$$

$$
\sum_{u=1}^{\infty}(-1)^{u}\frac{2u^2+1/2}{(2u^2-1/2)^2}=-1 \;\;\;\;\;\;\;\;\;\;\;\;\;\;\;\;\;\;\;\;\;\;
$$

$$
\sum_{u=1}^{\infty}\frac{(4^u-1)(u-1/2)-1}{2^{u^2+u+1}}=\frac{1}{4} \;\;\;\;\;\;\;\;\;\;\;\;\;\;\;\;\;\;\;\;
$$

$$
\sum_{u=1}^{\infty}\frac{(u^2+1/4)\tan(1/2)\cos u-u\sin u}{(4u^2-1)^2}=-\frac{\tan(1/2)}{8}
$$
with the generating functions, respectively
$$
F(n)=\frac{-1}{2(2n-1)}, \;\;\; F(n)=\frac{(-1)^{n-1}}{(2n-1)^2}, 
$$
$$
F(n)=\frac{-(n-1/2)}{2^{n^2-n+1}}, \;\;\; F(n)=\frac{\sin(n-1/2)}{8(2n-1)^2\cos(1/2)}
$$

\item[(2)] Divergent series
$$
1-1+1-1+...=\frac{1}{2} \;\;\;\;\;\;\;\;\;\;\;\;\;\;\;\;\;\;\;\;\;\;\;\;\;\; \textup{\citep{G}}
$$

\begin{equation}
1+1+1+1+...=-\frac{1}{2} \;\;\;\;\;\;\;\;\;\;\;\;\;\;\;\;\;\;\;\;\;\; \textup{\citep{H}} \label{eq:sum_of_ones}
\end{equation}

\begin{equation}
1^{2k}+2^{2k}+3^{2k}...=0 \;\;\;\;\; \forall k \;\;\;\;\;\;\;\;\;\;\;\;\;\;\;\;\;\;\;\; \textup{\citep{H}} \label{eq:sum_even_powers}
\end{equation}

\begin{equation}
1^{2k}-2^{2k}+3^{2k}-...=0 \;\;\;\;\; \forall k \;\;\;\;\;\;\;\;\;\;\;\;\;\;\;\;\;\;\; \textup{\citep{G}} \label{eq:sum_alt_even_powers}
\end{equation}

\begin{equation}
\sum_{u=1}^{\infty}(-1)^{u-1}u^{2k-1}\sin u\theta=0 \;\; \forall k, \;\;\; -\pi<\theta<\pi \;\; \textup{\citep{G}}  \label{ex:diver_sin}
\end{equation}

$$
\sum_{u=1}^{\infty}(-1)^{u-1}u^{2k}\cos u\theta=0 \;\;\;\;\; \forall k, \;\;\; -\pi<\theta<\pi \;\;\;\;\;\;\;\;\; \textup{\citep{G}}
$$
with the generating functions, respectively
$$
F(n)=\frac{(-1)^n}{2}, \;\;\;\;\; F(n)=n-1, \;\;\;\;\; F(n)=B_{2k}(n-1)
$$
$$
F(n)=\frac{(-1)^n}{2k+1}\sum_{u=1}^{2k+1}(2^u-1)\biggl(\begin{matrix}2k+1\\ u\end{matrix}\biggr)B_{u}(n-1)^{2k+1-u}
$$
$$
F(n)=(-1)^n\left(\sum_{u=1}^{k}\beta_{u}(n-1/2)^{2u-1}\sin(n-1/2)\theta +\sum_{u=0}^{k-1}\bar{\beta_{u}}(n-1/2)^{2u}\cos(n-1/2)\theta\right)
$$
$$
F(n)=(-1)^n\left(\sum_{u=1}^{k}\gamma_{u}(n-1/2)^{2u-1}\sin(n-1/2)\theta +\sum_{u=0}^{k}\bar{\gamma_{u}}(n-1/2)^{2u}\cos(n-1/2)\theta\right)
$$
where the coefficients $\beta_u, \bar{\beta_u}, \gamma_u, \bar{\gamma_u}$ are chosen such that the equalities hold
$$
F(n+1)-F(n)=(-1)^{n-1}n^{2k-1}\sin n\theta
$$
and
$$
F(n+1)-F(n)=(-1)^{n-1}n^{2k}\cos n\theta
$$
\end{itemize}

For $\theta=\frac{\pi}{2}$ in (\ref{ex:diver_sin}), we get
\begin{equation}
1^{2k-1}-3^{2k-1}+5^{2k-1}-...=0 \;\;\;\;\; \forall k \;\;\;\;\;\;\;\;\;\;\;\;\;\;\;\; \textup{\citep{G}}  \label{eq:from_diver_sin}
\end{equation}
\end{example}

\begin{corollary}[Theorem \ref{theorem:inf_sum}]
For any regular function $f(x)$ the following equality holds
\begin{equation}
\sum_{u=1}^{\infty}\left(f(u)+f(-u)\right)=-f(0)  \label{eq:1_cor_theor_inf_sum}
\end{equation}
independently on whether the series is convergent or not.
\end{corollary}

\begin{example}
It is easily checked that the function $f(x)=\frac{1}{9x^2-3x-2}$ is regular and 
$$
f(x)+f(-x)=\frac{18x^2-4}{81x^4-45x^2+4}
$$
Hence, the sum of convergent series
$$
\sum_{u=1}^{\infty}\frac{18u^2-4}{81u^4-45u^2+4}=\frac{1}{2}
$$
In particular, if we take $f(x)=e^x$, we get
$$
\sum_{u=1}^{\infty}(e^u+e^{-u})=-1
$$
and since $\cosh x=\frac{e^x+e^{-x}}{2}$ we finally obtain
$$
\sum_{u=1}^{\infty}\cosh u=-\frac{1}{2}.
$$
\end{example}

\begin{theorem} \label{theorem:theta1}
Let $\theta(x)$ be an elementary function such that 
\begin{equation}
\theta(-x)=-\theta(x-\epsilon), \label{eq:1_theor_theta1}
\end{equation} 
where $\epsilon=0,1,-1$.
Then
\begin{equation}
\lim_{n\rightarrow\infty}\theta(n+\delta)=0 \;\;\;\;\; (\delta=(1-\epsilon)/2) \label{eq:2_theor_theta1}
\end{equation}
\end{theorem}
\begin{proof}
Let $\theta(x+\delta)-\theta(x-1+\delta)=f(x)$. Then, according to axiom (\ref{axiom:5}), for any natural number $n$
$$
\theta(n+\delta)-\theta(\delta)=\sum_{u=1}^{n}f(u).
$$
The function $f(x)$ is regular, and moreover, it is even, $f(-x)=f(x)$. Indeed, $f(-x)=\theta\left(-(x-\delta)\right)-\theta\left(-(x+1-\delta)\right)=-\theta(x-\delta-\epsilon)+\theta(x+1-\delta-\epsilon)$. And since $-\delta-\epsilon=-1+\delta$, $f(-x)=-\theta(x-1+\delta)+\theta(x+\delta)=f(x)$. Therefore, passing to the limit and using (\ref{eq:inf_sum}), we get
$$
\lim_{n\rightarrow\infty}\theta(n+\delta)+\theta(\delta)=\lim_{n\rightarrow\infty}\sum_{u=1}^{n}f(u)=\sum_{u=1}^{\infty}f(u)=-\frac{f(0)}{2}
$$
and since $f(0)=\theta(\delta)-\theta(-1+\delta)=\theta(\delta)+\theta(\delta)=2\theta(\delta)$, we finally obtain
$$
\lim_{n\rightarrow\infty}\theta(n+\delta)=0,
$$
which completes the proof.
\end{proof}

\begin{remark} \label{rem:of_sequences}
Let $k_1, k_2,...,k_r$ be respectively $m_1,m_2,...,m_r$ - sequences of natural numbers and let
$$
Q_r(n)=\prod_{i=1}^{r}\Big((k_i+1)B_{k_i}(n)\Big)^{m_i}=\sum_{u=1}^{N_r}q_{u}n^{u},
$$
where $N_r=\sum_{i=1}^{r}(k_i+1)m_i$. Then, if $N_r\equiv1\;(\!\!\!\!\mod2)$ we, in accordance to (\ref{eq:2_theor_theta1}) and (\ref{eq:bernoulli}), get
\begin{equation}
\lim_{n\rightarrow\infty} Q_r(n)=0. \label{eq:rem_qu_zero}
\end{equation}
\end{remark}

\begin{remark} \label{rem:omega1}
Let $\omega(x)$ be an elementary function and let
\begin{equation}
\omega(-x)=-\omega(x-\epsilon t),  \label{eq:rem_omega1}
\end{equation}
where $t$ is a fixed natural number and $\epsilon=\pm1$.
Then, it is not hard to show that the expression
$$
H_{\omega}(x)=\sum_{u=\delta}^{t-1+\delta}\omega(x-\epsilon u)
$$
satisfies the condition
$$
H_{\omega}(-x)=-H_{\omega}(x-1)
$$
and in view of (\ref{eq:2_theor_theta1})
$$
\lim_{n\rightarrow\infty}H_{\omega}(n)=0.
$$
\end{remark}

\begin{example}
The polynomial $\theta(x)=x^3+\frac{9}{2}x^2+x-\frac{21}{4}$ satisfies the condition (\ref{eq:rem_omega1}) for $t=3, \epsilon=1$, and $H_{\theta}(x)=\sum_{u=0}^{2}\theta(x-u)=3x^3+\frac{9}{2}x^2-9x-\frac{21}{4}$. Indeed, 
$$
\theta(-x)=-x^3+\frac{9}{2}x^2-x-\frac{21}{4},
$$
$$
\theta(x-3)=(x^3-9x^2+27x-27)+(\frac{9}{2}x^2-27x+\frac{81}{2})+(x-3)-\frac{21}{4}=x^3-\frac{9}{2}x^2+x+\frac{21}{4},
$$
i. e.
$$
\theta(-x)=-\theta(x-3)
$$
and
\begin{eqnarray}
&&H_{\theta}(x)=\biggl[x^3+\frac{9}{2}x^2+x-\frac{21}{4}\biggr]+\biggl[\Bigl(x^3-3x^2+3x-1\Bigr)+\Bigl(\frac{9}{2}x^2-9x+\frac{9}{2}\Bigr)+ \nonumber \\
&&+
\Bigl(x-1\Bigr)-\frac{21}{4}\biggr]+\biggl[\Bigl(x^3-6x^2+12x-8\Bigr)+\Bigl(\frac{9}{2}x^2-18x+18\Bigr)+ \Bigl(x-2\Bigr)-\frac{21}{4}\biggr]= \nonumber \\
&&=
3x^3+\frac{9}{2}x^2-9x-\frac{21}{4}\nonumber.
\end{eqnarray}
The function $H_{\theta}(x)=3x^3+\frac{9}{2}x^2-9x-\frac{21}{4}$ satisfies the condition 
$$
H_{\theta}(-x)=-H_{\theta}(x-1).
$$ 
Indeed,
$$
H_{\theta}(-x)=-3x^3+\frac{9}{2}x^2+9x-\frac{21}{4}
$$
and
$$
H_{\theta}(x-1)=\Bigl(3x^3-9x^2+9x-3\Bigr)+\Bigl(\frac{9}{2}x^2-9x+\frac{9}{2}\Bigr)-\Bigl(9x-9\Bigr)-\frac{21}{4}=3x^3-\frac{9}{2}x^2-9x+\frac{21}{4},
$$
i. e.
$$
H_{\theta}(-x)=-H_{\theta}(x-1).
$$
Hence, in view of (\ref{eq:2_theor_theta1}),
$$
\lim_{n\rightarrow\infty}H_{\theta}(n)=0.
$$
Indeed (see Theorem \ref{eq:theor_polynom2}),
$$
\lim_{n\rightarrow\infty}H_{\theta}(n)=
\int\limits_{-1}^{0}\Bigl(3x^3+\frac{9}{2}x^2-9x-\frac{21}{4}\Bigr)dx=-\frac{3}{4}+\frac{3}{2}+\frac{9}{2}-\frac{21}{4}=0.
$$
\end{example}

\begin{proposition} \label{prop:odd_psi}
For any odd elementary function $\psi(x)$ the following equality holds
\begin{equation}
\lim_{n\rightarrow\infty}\sum_{u=\delta}^{t-1+\delta}\psi(n+\frac{\epsilon t}{2}-\epsilon u)=0.  \label{eq:prop_odd_psi}
\end{equation}
\end{proposition}
\begin{proof}
The function $\psi(x+\frac{\epsilon t}{2})$, where $\psi(x)$ is an odd elementary function, satisfies the condition (\ref{eq:rem_omega1}). Indeed,
$$
\psi(-x+\frac{\epsilon t}{2})=\psi\Bigl(-(x-\frac{\epsilon t}{2})\Bigr)=-\psi(x-\frac{\epsilon t}{2}).
$$
Therefore, using Remark \ref{rem:omega1} we obtain (\ref{eq:prop_odd_psi}).
\end{proof}

\begin{remark} \label{rem:odd_psi}
Relying on (3.12), we easily obtain
\begin{equation}
\lim_{n\rightarrow\infty}\psi(n+\frac{1}{2})=0,   \label{eq:1_rem_odd_psi}
\end{equation}
which is in good agreement with (2.17),
and also
$$
\lim_{n\rightarrow\infty}\psi(n+\frac{a}{2})=-\lim_{n\rightarrow\infty}\psi(n+1-\frac{a}{2})\;\;\;\;\;\forall a\in \mathbb{Z}
$$
or, in particular
\begin{equation}
\lim_{n\rightarrow\infty}(2n+1)^{2k-1}=0 \;\;\;\;\; \forall k\in \mathbb{N}  \label{eq:2_rem_odd_psi}
\end{equation}
and for $k=1$
\begin{equation}
\lim_{n\rightarrow\infty}n=-\frac{1}{2}  \label{eq:3_rem_odd_psi}
\end{equation}
\end{remark}

\begin{theorem} \label{theor:vartheta}
Let $\vartheta(x)$ be an elementary function such that
\begin{equation}
\vartheta(-x)=\vartheta(x-\epsilon t)\;\;\;\;\;\;\;(\epsilon=\pm 1)  \label{eq:theor_vartheta}
\end{equation}
Then the following equalities hold
$$
\lim_{n\rightarrow\infty}(-1)^{n}\vartheta(n+\delta)=0
$$
and
$$
\lim_{n\rightarrow\infty}(-1)^{n}\sum_{u=\delta}^{t-1+\delta}\vartheta(n-\epsilon u)=0.
$$
\end{theorem}
\begin{proof}
The proof is in the same way and uses the same techniques as for previous propositions.
\end{proof}

\begin{corollary}[Theorem \ref{theor:vartheta}] \label{cor:theor_vartheta}
For any even elementary function $\mu(x)$
\begin{equation}
\lim_{n\rightarrow\infty}(-1)^{n}\sum_{u=\delta}^{t-1+\delta}\mu(n+\frac{\epsilon t}{2}-\epsilon u)=0 \;\;\;\;\;\;\; (\epsilon=\pm 1). \label{eq:cor_theor_vartheta}
\end{equation}
\end{corollary}

\begin{remark} \label{rem:cor_vartheta}
Relying on (\ref{eq:cor_theor_vartheta}) we get
\begin{equation}
\lim_{n\rightarrow\infty}(-1)^{n}\mu(n+\frac{1}{2})=0,  \label{eq:1_rem_cor_vartheta}
\end{equation}
$$
\lim_{n\rightarrow\infty}(-1)^{n}\mu(n+\frac{a}{2})=-\lim_{n\rightarrow\infty}(-1)^{n}\mu(n+1-\frac{a}{2})\;\;\;\;\; \forall a\in \mathbb{Z}
$$
and
\begin{equation}
\lim_{n\rightarrow\infty}(-1)^{n}(2n+1)^{2k}=0 \;\;\;\;\; \forall k\in \mathbb{N}  \label{eq:2_rem_cor_vartheta}
\end{equation}
\end{remark}

\begin{example}
Let us take the well-known equality
$$
\sum_{u=1}^{n}\cos u\theta=-\frac{1}{2}+\Bigl(2\sin\frac{\theta}{2}\Bigr)^{-1}\sin\Bigl(n+\frac{1}{2}\Bigr)\theta,
$$
which holds for any natural $n$ and real $\theta$. Then, according to axiom (\ref{axiom:1}), we get
$$
\sum_{u=1}^{\infty}\cos u\theta=-\frac{1}{2}+\Bigl(2\sin\frac{\theta}{2}\Bigr)^{-1}\lim_{n\rightarrow\infty}\sin\Bigl(n+\frac{1}{2}\Bigr)\theta
$$
or, in view of (\ref{eq:1_rem_odd_psi}), 
\begin{equation}
\sum_{u=1}^{\infty}\cos u\theta=-\frac{1}{2}. \;\;\;\;\;\;\;\;\;\;\;\;\;\;\;\;\;\;\;\;\;\;\;\;\;\; \textup{\citep{G}} \label{ex:cos}
\end{equation}
Now taking the elementary trigonometric identity
$$
\sum_{u=1}^{n}(-1)^{u-1}\cos u\theta=\frac{1}{2}-\Bigl(2\cos\frac{\theta}{2}\Bigr)^{-1}(-1)^{n}\cos\Bigl(n+\frac{1}{2}\Bigr)\theta,
$$
passing to the limit and taking into account (\ref{eq:1_rem_cor_vartheta}), we obtain
\begin{equation}
\sum_{u=1}^{\infty}(-1)^{u-1}\cos u\theta=\frac{1}{2}. \;\;\;\;\;\;\;\;\;\;\;\;\;\;\;\;\;\;\;\;\;\; \textup{\citep{G}}  \label{ex:alt_cos}
\end{equation}
\end{example}
The last two formulas can also be obtained with use of (\ref{eq:inf_sum}).

\begin{theorem} \label{theor:gen_inf_sum}
For any regular function $f(x)$ satisfying the condition (\ref{eq:theor_vartheta}) the following equality holds
\begin{equation}
\sum_{u=1}^{\infty}f(u)=\frac{\epsilon}{2}\sum_{u=\delta}^{t-1+\delta}\Bigl(\lim_{n\rightarrow\infty}f(n-\epsilon u)-f(-\epsilon u)\Bigr)-\frac{1}{2}f(0). \label{eq:gen_inf_sum}
\end{equation}
\end{theorem}
\begin{proof}
According to formula (\ref{eq:1_from_lemma3}),
$$
\sum_{u=1}^{-(n+1)}f(u)=-\sum_{u=0}^{n}f(-u).
$$
In view of (\ref{eq:theor_vartheta}),
$$
\sum_{u=0}^{n}f(-u)=\sum_{u=0}^{n}f(u-\epsilon t)=\sum_{u=1}^{n}f(u)-\epsilon\sum_{u=\delta}^{t-1+\delta}f(n-\epsilon u)+\epsilon\sum_{u=\delta}^{t-\delta}f(-\epsilon u).
$$
Since 
$$
\sum_{u=\delta}^{t-\delta}f(-\epsilon u)=\sum_{u=\delta}^{t-1+\delta}f(-\epsilon u)+\epsilon f(-\epsilon t)
$$ 
and $f(-\epsilon t)=f(0)$, we have
$$
\sum_{u=0}^{n}f(-u)=\sum_{u=1}^{n}f(u)-\epsilon\sum_{u=\delta}^{t-1+\delta}\Bigl(f(n-\epsilon u)-f(-\epsilon u)\Bigr)+f(0).
$$
Hence,
$$
\sum_{u=1}^{-(n+1)}f(-u)=-\sum_{u=1}^{n}f(u)+\epsilon\sum_{u=\delta}^{t-1+\delta}\Bigl(f(n-\epsilon u)-f(-\epsilon u)\Bigr)-f(0).
$$
Passing to the limit and taking into account (\ref{eq:1_cor_lemma2}), we finally obtain
$$
\sum_{u=1}^{\infty}f(u)=\frac{\epsilon}{2}\sum_{u=\delta}^{t-1+\delta}\Bigl(\lim_{n\rightarrow\infty}f(n-\epsilon u)-f(-\epsilon u)\Bigr)-\frac{1}{2}f(0).
$$
The proof is completed.
\end{proof}
According to Theorem \ref{theor:gen_inf_sum}, for $t=1$ we have
$$
\sum_{u=1}^{\infty}f(u)=\frac{1}{2}\lim_{n\rightarrow\infty}f(n)-f(0) \;\;\;\;\; (\epsilon=1),
$$
\begin{equation}
\sum_{u=1}^{\infty}f(u)=-\frac{1}{2}\lim_{n\rightarrow\infty}f(n+1) \;\;\;\;\; (\epsilon=-1)  \label{eq:1_theor_gen_inf_sum}
\end{equation}
or, in particular, for $N_r\equiv 0 \:(\!\!\!\!\mod 2)$, since $Q(0)=0$, we have
\begin{equation}
\lim_{n\rightarrow\infty}Q_r(n)=2\sum_{u=1}^{N_r}q_u \sum_{\nu=1}^{\infty}\nu^u.  \label{eq:2_theor_gen_inf_sum}
\end{equation}

\begin{remark}
Formally, the formula (\ref{eq:gen_inf_sum}) coincides with (\ref{eq:inf_sum}) for $\epsilon=0$. So, the Theorem \ref{theor:gen_inf_sum} can be considered as a generalization of Theorem \ref{theorem:inf_sum} for the class of quasi-even functions, i. e. for the set of all regular functions $f(x)$ such that $f(-x)=f(x-a)$, where $a$ is a fixed integer number.
\end{remark}

Using the proof of the Theorem \ref{theor:gen_inf_sum}, we obtain the following
\begin{corollary}[Theorem \ref{theor:gen_inf_sum}] \label{cor:gen_theor}
For any regular function $f(x)$ the following equality holds
\begin{equation}
\sum_{u=1}^{\infty}f(u-\epsilon t)=\sum_{u=1}^{\infty}f(u)-\epsilon \sum_{u=\delta}^{t-1+\delta}\Bigl(\lim_{n\rightarrow\infty}f(n-\epsilon u)-f(-\epsilon u)\Bigr). \label{eq:cor_gen_theor}
\end{equation}
\end{corollary}

So, let $a_1(x)$, $a_2(x)$,..., $a_r(x)$ are regular functions, $\alpha_1$, $\alpha_2$,..., $\alpha_r$ are real numbers and $\sum_{i=1}^{r}\alpha_i \: a_i(x)=f(x-\epsilon t)$. Then
\begin{eqnarray}
\sum_{u=1}^{\infty}f(u)&=&\alpha_1 \sum_{u=1}^{\infty}a_1(u)+\alpha_2 \sum_{u=1}^{\infty}a_2(u)+...
+ \alpha_r \sum_{u=1}^{\infty}a_r(u)+ \nonumber \\
&&+\;
\epsilon \sum_{u=\delta}^{t-1+\delta}\Bigl(\lim_{n\rightarrow\infty}f(n-\epsilon u)-f(-\epsilon u)\Bigr) \nonumber.
\end{eqnarray}

\begin{example}
Thus, letting $a_1(x)=x^2$, $a_2(x)=x$ and $a_3(x)=1$, we have
$a_1(x)+2a_2(x)+a_3(x)=(x+1)^2=f(x+1)$, where $f(x)=x^2$, $\epsilon t=-1$. Then
$$
\sum_{u=1}^{\infty}u^2=\sum_{u=1}^{\infty}u^2+2\sum_{u=1}^{\infty}u+\sum_{u=1}^{\infty}1-\lim_{n\rightarrow\infty}(n+1)^2 +1
$$
or, taking into account (\ref{eq:sum_of_1}) and (\ref{eq:3_rem_odd_psi})
$$
\sum_{u=1}^{\infty}u=\frac{1}{2}\lim_{n\rightarrow\infty}n^2+\lim_{n\rightarrow\infty}n+\frac{1}{4}
$$
and
\begin{equation}
\sum_{u=1}^{\infty}u=\frac{1}{2}\lim_{n\rightarrow\infty}n^2-\frac{1}{4}. \label{ex:cor_gen_theor}
\end{equation}
\end{example}

\begin{remark} \label{rem:gen_theor}
If \: $\sum_{i=1}^{r}\alpha_i \: a_i(x)=f(x)$, i. e. $\epsilon=0$, then
\begin{equation}
\sum_{u=1}^{\infty}f(u)=\alpha_1 \sum_{u=1}^{\infty}a_1(u)+\alpha_2 \sum_{u=1}^{\infty}a_2(u)+...+ \alpha_r \sum_{u=1}^{\infty}a_r(u) \label{eq:rem_gen_theor}
\end{equation}
and since $4a_1(x)-4a_2(x)+a_3(x)=(2x-1)^2=f(x)$, $f(x)=(2x-1)^2$, $\epsilon=0$
$$
\sum_{u=1}^{\infty}(2u-1)^2=4\sum_{u=1}^{\infty}u^2-4\sum_{u=1}^{\infty}u+\sum_{u=1}^{\infty}1
$$
and, according to (\ref{eq:sum_of_1}) and (\ref{eq:cor_gen_theor})
$$
\sum_{u=1}^{\infty}(2u-1)^2=\frac{1}{2}-2\lim_{n\rightarrow\infty}n^2.
$$
Analogously, we get
$$
\sum_{u=1}^{\infty}(2u-1)=2\sum_{u=1}^{\infty}u-\sum_{u=1}^{\infty}1
$$
and, according to (\ref{eq:sum_of_1}) and (\ref{eq:cor_gen_theor})
\begin{equation}
\sum_{u=1}^{\infty}(2u-1)=\lim_{n\rightarrow\infty}n^2.  \label{eq:rem_gen_theor_lim}
\end{equation}
\end{remark}
The functions $f(n)=(2n-1)^{2k}$ and $f(n)=(-1)^{n}(2n-1)^{2k-1}$, \; ($\epsilon t=-1$), satisfies the condition of the Theorem \ref{theor:gen_inf_sum}. Therefore, in view of (\ref{eq:1_theor_gen_inf_sum})
\begin{equation}
\sum_{u=1}^{\infty}(2u-1)^{2k}=-\frac{1}{2}\lim_{n\rightarrow\infty}(2n+1)^{2k} \;\;\;\;\;\;\; \forall k   \label{eq:1_from_distrib}
\end{equation}
$$
\sum_{u=1}^{\infty}(-1)^{u-1}(2u-1)^{2k-1}=-\frac{1}{2}\lim_{n\rightarrow\infty}(-1)^{n}(2n+1)^{2k-1} \;\;\;\;\;\;\; \forall k
$$
and, according to (\ref{eq:from_diver_sin})
\begin{equation}
\lim_{n\rightarrow\infty}(-1)^{n}(2n+1)^{2k-1}=0 \;\;\;\;\;\;\; \forall k.  \label{eq:2_from_distrib}
\end{equation}

Combining the formulas (\ref{eq:2_rem_cor_vartheta}) and (\ref{eq:2_from_distrib}), we arrive at the following
\begin{proposition}
For any natural number $k$
\begin{equation}
\lim_{n\rightarrow\infty}(-1)^{n}(2n+1)^{k}=0.  \label{eq:alt_odd_int_zero}
\end{equation}
\end{proposition}

\subsection{Theorems on Polynomials}

\begin{theorem} \label{theor:polynom1}
Suppose $f(x)$ is a polynomial defined over the field of real numbers, $x\in \mathbb{R}$. Then, the following equality holds
\begin{equation}
\lim_{n\rightarrow\infty}(-1)^{n}f(n)=0.  \label{eq:theor_polynom1}
\end{equation}
\end{theorem}
\begin{proof}
To prove the theorem,  it is sufficient to show, due to the rule \ref{rule2}, 
that for any non-negative integer value $\sigma$
\begin{equation}
\lim_{n\rightarrow\infty}(-1)^{n}n^{\sigma}=0.  \label{sigma}
\end{equation}
The proof is by induction over $\sigma$. For $\sigma=0$ the formula (\ref{sigma}), according to (\ref{eq:1_rem_cor_vartheta}), holds true. Assume now that it holds for all positive integers less than some natural number $k$. Then, using (\ref{eq:alt_odd_int_zero}), we get the expansion
$$
2^k\lim_{n\rightarrow\infty}(-1)^{n}n^k+2^{k-1}\biggl(\begin{matrix}k\\ 1\end{matrix}\biggr)\lim_{n\rightarrow\infty}(-1)^{n}n^{k-1}+...+\lim_{n\rightarrow\infty}(-1)^{n}=0,
$$
in which, by the inductive assumption, all the terms, except for the first, are equal to zero. Hence, the first term will be equal to zero as well, i. e.
$$
2^k\lim_{n\rightarrow\infty}(-1)^{n}n^k=0
$$
and
$$
\lim_{n\rightarrow\infty}(-1)^{n}n^k=0.
$$
This means, by virtue of the induction, that the formula (\ref{sigma}) holds for all non-negative integer values of $\sigma$. 

The theorem is proved.
\end{proof}

The following theorem demonstrates a connection between the limit of polynomial with the integral of it.

\begin{theorem} \label{theor:polynom2}
Suppose $f(x)$ is a polynomial defined over the field of real numbers, $x\in \mathbb{R}$. Then, the following equality holds
\begin{equation}
\lim_{n\rightarrow\infty}f(n)=\int\limits_{-1}^{0}f(x)dx. \label{eq:theor_polynom2}
\end{equation}
\end{theorem}
\begin{proof}
To prove the theorem,  it is sufficient to show, analogously to the previous theorem, that the relation
\begin{equation}
\lim_{n\rightarrow\infty}n^{\sigma}=\frac{(-1)^{\sigma}}{\sigma+1} \label{sigma_fraction}
\end{equation}
holds for any non-negative integer value $\sigma$.

We give here a sketch of the proof.

At first step, we have the following equality
$$
\sum_{u=1}^{n}u^{2k}+\sum_{u=1}^{n}(-1)^{u-1}u^{2k}=2\sum_{u=1}^{[n/2]}(2u-1)^{2k}+\Bigl(n^{2k}-(-1)^{n}n^{2k}\Bigr) \;\;\;\;\;\;\; \forall \: n, k.
$$
Passing to the limit and taking into account the axioms (\ref{axiom:1}), (\ref{axiom:2}) and formula (\ref{sigma}), we obtain the relation
$$
\sum_{u=1}^{\infty}u^{2k}+\sum_{u=1}^{\infty}(-1)^{u-1}u^{2k}=2\sum_{u=1}^{\infty}(2u-1)^{2k}+\lim_{n\rightarrow\infty}n^{2k}
$$
or, according to (\ref{eq:sum_even_powers}), (\ref{eq:sum_alt_even_powers}) and (\ref{eq:1_from_distrib})
\begin{equation}
\lim_{n\rightarrow\infty}(2n+1)^{2k}=\lim_{n\rightarrow\infty}n^{2k}.  \label{eq:lim_even_powers}
\end{equation}
Then, it is proved that 
\begin{equation}
\sum_{u=0}^{2k}(-1)^{u}\: 2^{2k-u}\biggl(\begin{matrix}2k\\ u\end{matrix}\biggr)\frac{1}{2k+1-u}=\frac{1}{2k+1} \label{eq:1_proof_polynom2}
\end{equation}
and
\begin{equation}
\sum_{u=0}^{2k-1}(-1)^{u}\: 2^{2k-1-u}\biggl(\begin{matrix}2k-1\\ u\end{matrix}\biggr)\frac{1}{2k-u}=0. \label{eq:2_proof_polynom2}
\end{equation}
Then, following the induction over $\sigma$, we assume that the equality (\ref{sigma_fraction}) holds for all non-negative integers less than some natural number $k$. Further, depending on the parity of $k$ we use the formulas (\ref{eq:1_proof_polynom2}) and (\ref{eq:lim_even_powers}) or formulas (\ref{eq:2_proof_polynom2}3.42) and (\ref{eq:2_rem_odd_psi}), and then prove that (\ref{sigma_fraction}) holds for $\sigma=k$. And since the formula (\ref{sigma_fraction}) holds for $\sigma=0$, by the induction hypothesis it holds for all non-negative integer values of $\sigma$.

This completes the proof.
\end{proof}

Below we reprove the known theorem on the binomial series and obtain new result, proving the validity of the theorem at the right endpoint
of the circle of convergence.

\begin{theorem} \label{theor:binom}
For any integer $a\neq 0$ and real $x$, $-1<x\leq 1$ the following equality holds
\begin{equation}
(1+x)^a = \sum_{u=0}^{\infty}x^{u}\:\frac{\lambda(a)}{\lambda(u)\lambda(a-u)}= 1+x\biggl(\:\begin{matrix}a\\ 1\end{matrix}\:\biggr)+x^{2}\biggl(\:\begin{matrix}a\\ 2\end{matrix}\:\biggr)+... \label{eq:theor_binom}
\end{equation}
\end{theorem}
\begin{proof}
For $a>0$ the formula (\ref{eq:theor_binom}) is reduced to the elementary binomial formula. Thus, let $a$ be a negative number, that is $a=-m$. In this case, dividing $1$ by $(1+x)^m$, using the rule of division of one polynomial by another, we will gradually obtain in the quotient the terms of the right-hand side of (\ref{eq:theor_binom}). 
Stopping the process of division at some fixed number $k$, we get
$$
(1+x)^{-m} = 1-x\biggl(\:\begin{matrix}m\\ 1\end{matrix}\:\biggr)+x^{2}\biggl(\:\begin{matrix}m-1\\ 2\end{matrix}\:\biggr)-...(-1)^{k}x^{k}\biggl(\:\begin{matrix}m-k-1\\ k\end{matrix}\:\biggr)+R_{k}^{m}(x),
$$
where
$$
R_{k}^{m}(x)=\frac{(-1)^{k+1}}{(1+x)^m}\sum_{u=0}^{m-1}\biggl(\:\begin{matrix}k+u\\ u\end{matrix}\:\biggr)\biggl(\:\begin{matrix}m+k\\ m-1-u\end{matrix}\:\biggr)x^{k+1+u}.
$$
If $|x|<1$ and $k$ is large enough, then obviously $R_{k}^{m}(x)$ is arbitrarily small and it tends to zero as $k$ unboundedly increases, that is
$$
\lim_{k\rightarrow\infty}R_{k}^{m}(x)=0 \;\;\;\;\;\;\; |x|<1.
$$
If $x=1$, then $\bigl|R_{k}^{m}(x)\bigr|$ becomes a polynomial of degree $m-1$ in the variable $k$, which according to (\ref{eq:theor_polynom1}) also satisfy the equality
$$
\lim_{k\rightarrow\infty}R_{k}^{m}(1)=0.
$$
This means that (\ref{eq:theor_binom}) holds for all $x$, $-1<x\leq 1$.

The proof is completed.
\end{proof}

\begin{example}
\begin{eqnarray}
	(1+1)^{-1}=1-1+1-...=\frac{1}{2}  \nonumber \;\;\;\;\;\;\;\;\;\;\;\;\;\;\;\;\;\;\;\;\;\;\;\;\;\; \textup{\citep{G}} \\
	(1+1)^{-2}=1-2+3-...=\frac{1}{4} \nonumber \;\;\;\;\;\;\;\;\;\;\;\;\;\;\;\;\;\;\;\;\;\;\;\;\;\; \textup{\citep{G}} \\
	(1+1)^{-3}=1-3+6-10+...=\frac{1}{8}  \nonumber \;\;\;\;\;\;\;\;\;\;\;\;\;\;\;\;\;\;\;\;\;\;\;\;\;\; \textup{\citep{G}} 
\end{eqnarray}
and so on.
\end{example}

\begin{remark}
The validity of the Theorem \ref{theor:binom} can be also established, relying on the formulas (\ref{eq:combi}) and (\ref{eq:theor_polynom1}).
\end{remark}

\subsection{Several Properties of Divergent Series}

\begin{theorem}[commutative property] \label{theor:commut}
For any regular function $f(x)$
\begin{equation}
\sum_{u=1}^{\infty}f(u)=\sum_{u=1}^{\infty}f(2u-1)+\sum_{u=1}^{\infty}f(2u)+\frac{1}{2}\lim_{n\rightarrow\infty}\Bigl(f(n)-(-1)^nf(n)\Bigr). \label{eq:theor_commut}
\end{equation}
\end{theorem}
\begin{proof}
Let us take the equality
$$
\sum_{u=1}^{n}f(u)=\sum_{u=1}^{[n/2]}f(2u-1)+\sum_{u=1}^{[n/2]}f(2u)+\frac{1}{2}\Bigl(f(n)-(-1)^nf(n)\Bigr)\;\;\;\;\;\;\; \forall n
$$
Now, passing to the limit and taking into account axioms (\ref{axiom:1}) and (\ref{axiom:2}) we arrive at (\ref{eq:theor_commut}). This completes the proof.
\end{proof}

\begin{corollary}[Theorem \ref{theor:commut}] \label{cor:theor_commut}
If \: $\lim_{n\rightarrow\infty}(-1)^{n}f(n)=0$, then
$$
\sum_{u=1}^{\infty}f(u)=\sum_{u=1}^{\infty}f(2u-1)+\sum_{u=1}^{\infty}f(2u)+\frac{1}{2}\lim_{n\rightarrow\infty}f(n).
$$
\end{corollary}

\begin{theorem}[associative property] \label{theor:assoc}
For any regular function $f(x)$
\begin{equation}
\sum_{u=1}^{\infty}f(u)=\sum_{u=1}^{\infty}\Bigl(f(2u-1)+f(2u)\Bigr)+\frac{1}{2}\lim_{n\rightarrow\infty}\Bigl(f(n)+(-1)^{n}f(n)\Bigr). \label{eq:theor_assoc}
\end{equation}
\end{theorem}
\begin{proof}
The proof is analogous to that of the Theorem \ref{theor:commut}.
\end{proof}
\begin{example}
Consider two simple examples.
\begin{itemize}
\item[(1)]
$$
\sum_{u=1}^{\infty}1=\sum_{u=1}^{\infty}(1+1)+\frac{1}{2}\lim_{n\rightarrow\infty}\Bigl(1+(-1)^{n}\Bigr)
$$
and
\item[(2)]
$$
\sum_{u=1}^{\infty}(-1)^{u-1}=\sum_{u=1}^{\infty}(1-1)+\frac{1}{2}\lim_{n\rightarrow\infty}\Bigl((-1)^{n}+1\Bigr).
$$
Whence, in view of (\ref{eq:1_rem_cor_vartheta}), we find
$$
\sum_{u=1}^{\infty}1=2\sum_{u=1}^{\infty}1+\frac{1}{2},
$$
i. e.
$$
\sum_{u=1}^{\infty}1=-\frac{1}{2}
$$
and 
$$
\sum_{u=1}^{\infty}(-1)^{u-1}=0+\frac{1}{2}=\frac{1}{2}.
$$
\end{itemize}
The last two equalities are in good agreement with the formula (\ref{eq:inf_sum}).
\end{example}

\begin{theorem} \label{theor:distrib}
Let $f(x)$ be a regular function and let $\sum_{u=1}^{\infty}f(u)=A$ and $\sum_{u=1}^{\infty}(-1)^{u-1}f(u)=B$. Then
\begin{equation}
\sum_{u=1}^{\infty}f(2u)=\frac{1}{2}(A-B)  \label{eq:1_theor_distrib}
\end{equation}
and
\begin{equation}
\sum_{u=1}^{\infty}f(2u-1)=\frac{1}{2}(A+B)-\frac{1}{2}\lim_{n\rightarrow\infty}f(n). \label{eq:2_theor_distrib}
\end{equation}
\end{theorem}
\begin{proof}
The proof uses the same techniques as in previous theorems.
\end{proof}



\subsection{Arithmetic and Geometric Progressions}

\begin{theorem} \label{theor:arith}
Let $a_1$, $a_2$,...,$a_u$,... be an infinite arithmetic progression, i. e. $a_u=a_1+(u-1)d$ and $d\geq 0$. Then for the arithmetic series the following equality holds
\begin{equation}
\sum_{u=1}^{\infty}a_u=\frac{5d-6a_1}{12}. \label{eq:theor_arith}
\end{equation}
\end{theorem}
\begin{proof}
The sum $S_n$ of the first $n$ terms of arithmetic progression with difference $d$ and first term $a_1$ is defined by the formula
$$
S_n=\sum_{u=1}^{n}a_u=\frac{2a_1+(n-1)d}{2}\:n.
$$
Then, passing to the limit and using (\ref{eq:theor_polynom2}), we get
$$
\lim_{n\rightarrow\infty}S_n = \lim_{n\rightarrow\infty}\sum_{u=1}^{n}a_u = \int\limits_{-1}^{0}\biggl(\frac{2a_1+(x-1)d}{2}\biggr)xdx = \frac{5d-6a_1}{12}
$$
and, in view of axiom (\ref{axiom:1}), we finally obtain
$$
\lim_{n\rightarrow\infty}S_n = \sum_{u=1}^{\infty}a_u=\frac{5d-6a_1}{12}.
$$
The theorem is proved.
\end{proof}

\begin{example}
\begin{eqnarray}
	\sum_{u=1}^{\infty}1=1+1+1+...=-\frac{1}{2}, \;\;\;\;\;\;\; (d=0)  \nonumber \;\;\;\;\;\;\;\;\;\;\;\;\;\;\;\;\;\;\;\;\;\;\;\;\;\; \textup{\citep{H}} \\ 
	\sum_{u=1}^{\infty}u=1+2+3+...=-\frac{1}{12},\;\;\;\;\;\;\; (d=1)  \nonumber \;\;\;\;\;\;\;\;\;\;\;\;\;\;\;\;\;\;\;\;\;\;\;\;\;\; \textup{\citep{H}} \\
	\sum_{u=1}^{\infty}(2u-1)=1+3+5+...=\frac{1}{3},\;\;\;\;\;\;\; (d=2) \nonumber \;\;\;\;\;\;\;\;\;\;\;\;\;\;\;\;\;\;\;\;\;\;\;\;\;\; \textup{\citep{G}}
\end{eqnarray}
and so on.
\end{example}

From Theorem \ref{theor:polynom1} we immediately obtain
\begin{theorem} \label{theor:alphabeta}
Let $\alpha(x)$ and $\beta(x)$ be elementary functions defined on $\mathbb{Z}$ and satisfying the condition $\alpha(x)-\beta(x)=f(x)$, where $f(x)$ is a polynomial. Suppose that $\mu(x)$ is a function such that $\mu(x)=\alpha(x)$ if \, $2\mid x$ and $\mu(x)=\beta(x)$ if \, $2/\!\!\!\!\mid x$. Then
$$
\lim_{n\rightarrow\infty}\mu(n)=\frac{1}{2}\lim_{n\rightarrow\infty}\bigl(\alpha(n)+\beta(n)\bigr).
$$
\end{theorem}
\begin{proof}
It is readily seen that function $\mu(x)$ satisfying the conditions of the theorem can be represented as the sum of two functions
$$
\mu(n)=\mu_{1}(n)+(-1)^{n} \: \mu_{2}(n),
$$
where $\mu_{1}(n)=\frac{1}{2}\bigl(\alpha(n)+\beta(n)\bigr)$ and
$\mu_{2}(n)=\frac{1}{2}\bigl(\alpha(n)-\beta(n)\bigr)$.\\
Therefore,
$$
\lim_{n\rightarrow\infty}\mu(n)=\lim_{n\rightarrow\infty}\bigl(\mu_{1}(n)+(-1)^{n}\:\mu_{2}(n)\bigr).
$$
But, by the circumstances of the theorem,
$$
\mu_{2}(n)=\frac{1}{2}\bigl(\alpha(n)-\beta(n)\bigr)=\frac{1}{2}f(n).
$$
Hence
$$
\lim_{n\rightarrow\infty}\mu(n)=\frac{1}{2}\lim_{n\rightarrow\infty}\bigl(\alpha(n)+\beta(n)\bigr)+\frac{1}{2}\lim_{n\rightarrow\infty}(-1)^{n}f(n)
$$
and, in view of (\ref{eq:theor_polynom1})
$$
\lim_{n\rightarrow\infty}\mu(n)=\frac{1}{2}\lim_{n\rightarrow\infty}\bigl(\alpha(n)+\beta(n)\bigr).
$$
The proof is completed.
\end{proof}

\begin{remark} 
As one can see from the proof of Theorem \ref{theor:alphabeta}, it actually holds for a considerably wider class of elementary functions $\alpha(x)$ and $\beta(x)$ which satisfy the condition
$$
\lim_{n\rightarrow\infty}(-1)^{n}\bigl(\alpha(n)-\beta(n)\bigr)=0.
$$
\end{remark}

Using Theorem \ref{theor:alphabeta}, we obtain the following

\begin{theorem} \label{theor:alt_arith}
Let $a_1$, $a_2$,...,$a_u$,... be an infinite arithmetic progression, i. e. $a_u=a_1+(u-1)d$ and $d\geq 0$. Then for the alternating arithmetic series the following equality holds
\begin{equation}
\sum_{u=1}^{\infty}(-1)^{u-1}\: a_u=\frac{2a_1-d}{4}. \label{eq:theor_alt_arith}
\end{equation}
\end{theorem}
\begin{proof}
Let
$$
\overline{S}_n=\sum_{u=1}^{n}(-1)^{u-1} \: a_u \;\;\;\;\;\;\; \forall \: n.
$$
If $n=2n_1$, then
$$
\overline{S}_n=\sum_{u=1}^{n_1} a^{\prime}_u - \sum_{u=1}^{n_1} b^{\prime}_u,
$$
where $a^{\prime}_u=a_1+(u-1)2d$ and $b^{\prime}_u=a_1+d+(u-1)2d$.\\
Therefore, 
$$
\overline{S}_n=\biggl(\frac{2a_1+(n_1-1)2d}{2}\biggr)n_1 - \biggl(\frac{2a_1+2d+(n_1-1)2d}{2}\biggr)n_1 = -d n_1
$$
and
$$
\overline{S}_n=-d \: \frac{n}{2} \;\;\;\;\;\;\;\;\;\;\;\; (2\mid x).
$$
If $n=2n_1+1$, then
$$
\overline{S}_n=\sum_{u=1}^{n_1+1}a^{\prime}_u - \sum_{u=1}^{n_1}b^{\prime}_u = a_1+dn_1
$$
and
$$
\overline{S}_n=a_1 + d \: \frac{n-1}{2} \;\;\;\;\;\;\;\;\;\;\;\; (2/\!\!\!\!\mid x).
$$
So, the function $\overline{S}_n$ takes the values of $-d(n/2)$ if $2\!\mid \!x$ and of $a_1+d(n-1)/2$ if $2/\!\!\!\!\mid x$. 
Hence, in view of Theorem \ref{theor:alphabeta}, since
$$
\biggl(-d\:\frac{n}{2}\biggr)-\biggl(a_1+d\:\frac{n-1}{2}\biggr)=-dn-\frac{2a_1-d}{2}  
$$
is a polynomial of $n$, we have
$$
\lim_{n\rightarrow\infty}\overline{S}_n=\lim_{n\rightarrow\infty}\sum_{u=1}^{n}(-1)^{u-1}a_u=\frac{1}{2}\lim_{n\rightarrow\infty}\biggl(a_1+\frac{d(n-1)}{2}-\frac{dn}{2}\biggr)=\frac{2a_1-d}{4}
$$
and finally
$$
\sum_{u=1}^{\infty}(-1)^{u-1}a_u=\frac{2a_1-d}{4}.
$$
The theorem is proved.
\end{proof}

\begin{example}
\begin{eqnarray}
	\sum_{u=1}^{\infty}(-1)^{u-1}=1-1+1-1+...=\frac{1}{2},\;\;\;\;\;\;\; (d=0)  \nonumber \;\;\;\;\;\;\;\;\;\;\;\;\;\;\;\;\;\;\;\;\; \textup{\citep{G}} \\
	\sum_{u=1}^{\infty}(-1)^{u-1}u=1-2+3-4+...=\frac{1}{4},\;\;\;\;\;\;\; (d=1) \nonumber \;\;\;\;\;\;\;\;\;\;\;\;\;\;\;\;\;\;\;\; \textup{\citep{G}} \\
	\sum_{u=1}^{\infty}(-1)^{u-1}(2u-1)=1-3+5-7+...=0,\;\;\;\;\;\;\; (d=2) \nonumber \;\;\;\;\;\;\;\;\;\;\;\;\;\;\;\;\;\;\;\; \textup{\citep{G}}
\end{eqnarray}
and so on.

The last equality is a special case, $k=1$, of the formula (\ref{eq:from_diver_sin}).
\end{example}

\begin{remark}
The Theorem \ref{theor:alt_arith} can be proved in another way by relying on axiom (\ref{axiom:2}). Indeed, let us take the equality
$$
\overline{S}_n + S_n = 2\sum_{u=1}^{[n/2]}a_{2u-1}+\bigl(1-(-1)^{n}\bigr)a_n,
$$
where $a_{2u-1}=a_1+(u-1)2d$, $a_n=a_1+(n-1)d$. Then, passing to the limit and taking into account (\ref{eq:3_rem_odd_psi}), (\ref{sigma}), 
Theorem \ref{theor:arith}, and axiom (\ref{axiom:2}), we get
$$
\lim_{n\rightarrow\infty}\overline{S}_n + \lim_{n\rightarrow\infty}S_n = 2\sum_{u=1}^{\infty}a_{2u-1}+\lim_{n\rightarrow\infty}a_n,
$$
i. e.
$$
\sum_{u=1}^{\infty}(-1)^{u-1}a_u + \frac{5d-6a_1}{12} = 2\biggl(\frac{10d-6a_1}{12}\biggr)+a_1-\frac{3}{2}d
$$
and
$$
\sum_{u=1}^{\infty}(-1)^{u-1}a_u = \frac{2a_1-d}{4}.
$$
\end{remark}

\begin{theorem} \label{theor:geom}
For any real number $g$, $g\neq 1$
\begin{equation}
\sum_{u=0}^{\infty}g^u = \frac{1}{1-g}. \label{eq:theor_geom}
\end{equation}
\end{theorem}
\begin{proof}
For $|g|<1$ the validity of the theorem is obvious. Thus, let $|g|>1$. Then $\bigl|\frac{1}{g}\bigr|<1$ and
$$
\sum_{u=1}^{\infty}\frac{1}{g^u}=\frac{1}{g-1}.
$$
It is clear that the function $\frac{1}{g^x}$ is regular. Hence, in view of (\ref{eq:2_from_lemma3})
$$
\sum_{u=1}^{\infty}\frac{1}{g^u}=-\sum_{u=0}^{\infty}g^u
$$
and
$$
\sum_{u=0}^{\infty}g^u=\frac{1}{1-g},
$$
which completes the proof.
\end{proof}

\begin{example}
\begin{eqnarray}
	\sum_{u=0}^{\infty}2^u = 1+2+4+...=-1  \nonumber \;\;\;\;\;\;\;\;\;\;\;\;\;\;\;\;\;\;\;\;\;\;\;\;\;\; \textup{\citep{G}} \\
	\sum_{u=0}^{\infty}(-1)^{u}\:2^u = 1-2+4-...=\frac{1}{3}  \nonumber \;\;\;\;\;\;\;\;\;\;\;\;\;\;\;\;\;\;\;\;\;\;\;\;\;\; \textup{\citep{G}} 
\end{eqnarray}
\end{example}

\vspace{1cm}
\begin{center}
\rule{9.5cm}{0.25mm}
\end{center}

\newpage

\section{Bernoulli Numbers and Some Identities}\label{sec:bern}

Bernoulli numbers play an important role in many areas of analysis and number theory. They often occur when expanding some simple functions in a power series. These numbers were first introduced by Jacob Bernoulli when he worked on the calculation of sums $1^k+2^k+3^k+...+n^k$. The Bernoulli nunbers can be defined by the generating series
$$
\frac{te^{t}}{e^{t}-1}=\frac{t}{1-e^{-t}}=\sum_{k=0}^{\infty}B_k\frac{t^k}{k!}.
$$
First several values of the Bernoulli numbers are \\

\begin{tabbing}

22222 \= 222222222222 \= 222222222222222 \= 222222222222 \= 222222222222 \= 222222222222  \kill

\> $B_0=1$, \> $B_1=\frac{1}{2}$, \> $B_2=\frac{1}{6}$, \> $B_3=0$, \> $B_4=-\frac{1}{30}$, \\ \\
\> $B_5=0$,  \> $B_6=\frac{1}{42}$, \> $B_7=0$,  \> $B_8=-\frac{1}{30}$, \>  $B_9=0$, \\ \\
\> $B_{10}=\frac{5}{66}$,  \> $B_{11}=0$, \> $B_{12}=\frac{691}{2730}$, \> $B_{13}=0$, \> $B_{14}=\frac{7}{6}$, \\ \\
\> $B_{15}=0$, \> $B_{16}=-\frac{3617}{510}$, \> $B_{17}=0$, \> $B_{18}=\frac{43867}{798}$, \> $B_{19}=0$, \\ \\ 
\> $B_{20}=-\frac{174611}{330}$,...  \\

\end{tabbing} 
Let us consider the Bernoulli polynomials, which can be defined by the symbolic equality $B_k(t)=(t+B)^k$, where after the brackets removing the exponents of $B$ should be substituted by indices, and $B_k$ are the Bernoulli numbers. Below we derive two equalities, well-known from analysis, without using the notions of functions of complex variable and analytic continuation.

On the one hand, according to (\ref{eq:bernoulli}) we have
$$
\frac{1}{k}B_k(-1)=\frac{-1}{k(k+1)}\sum_{u=1}^{k}(-1)^{u+k}\biggl(\begin{matrix}k+1\\ u\end{matrix}\biggr)B_u=\frac{(-1)^{k-1}}{k}\sum_{u=0}^{0}u^k=0.
$$

On the other hand, using (\ref{eq:theor_polynom2}) we get
$$
\lim_{n\rightarrow\infty}B_{k-1}(n)=\int\limits_{-1}^{0}B_{k-1}(x)dx=\frac{1}{k(k+1)}\sum_{u=0}^{k-1}(-1)^{u+k}\biggl(\begin{matrix}k+1\\ u\end{matrix}\biggr)B_u=\sum_{u=1}^{\infty}u^{k-1}.
$$

Adding these two expressions, we obtain
\begin{equation} \label{eq:analysis1}
\sum_{u=1}^{\infty}u^{k-1}=-\frac{B_k}{k} \;\;\;\;\;\;\;\;\; \forall k \;\;\;\;\;\;\;\;\;\;\;\;\; \textup{\citep{H}} 
\end{equation}     

Now taking the formula
$$
\frac{(2^k-1)B_k}{k}-\frac{(-1)^n}{k}\sum_{u=1}^{k}(2^u-1)\biggl(\:\begin{matrix}k\\ u\end{matrix}\:\biggr)B_{u} \: n^{k-u}=\sum_{u=1}^{n}(-1)^{u-1}u^{k-1} \;\;\;\;\; \forall k,
$$
which holds for any natural number $n$, and passing to the limit with regard to (\ref{eq:theor_polynom1}), we get
\begin{equation} \label{eq:analysis2}
\sum_{u=1}^{\infty}(-1)^{u-1}u^{k-1}=\frac{(2^k-1)B_k}{k} \;\;\;\;\;\;\;\;\; \forall k \;\;\;\;\;\;\;\;\; \textup{\citep{G}} 
\end{equation}

\subsection{Recurrence Formulas of Bernoulli Numbers}
   
Using the formulas (\ref{eq:gen_inf_sum}), (\ref{eq:theor_polynom1}), (\ref{eq:theor_polynom2}), (\ref{eq:analysis1}), (\ref{eq:analysis2}), and assuming $f(x)-f(x-1)=\sum_{u=0}^{k-1}a_u x^u$ and $f(x)+f(x-1)=\sum_{u=0}^{k}b_u x^u$, 
where $f(x)$ is a polynomial of degree $k$, we derive the following recurrence formulas for Bernoulli numbers:
$$
\sum_{u=0}^{k-1}\frac{a_u}{u+1}\:B_{u+1}=\int\limits_{0}^{-1}\bigl(f(x)-f(0)\bigr)dx,
$$
$$
\sum_{u=0}^{k}\frac{\bigl(2^{u+1}-1\bigr)b_u}{u+1}\:B_{u+1}=f(0),
$$
$$
\sum_{u=1}^{k}2^{2u-1}\bigl(2^{2u}-1\bigr)\biggl(\begin{matrix}2k\\2u\end{matrix}\biggr)B_{2u}=
\bigl(2^{2k}-1\bigr)B_{2k}+\sum_{u=0}^{k-1}2^{2u-1}\biggl(\begin{matrix}2k\\2u\end{matrix}\biggr)B_{2u}=k \;\;\;\;\; \forall k,
$$
\begin{equation}
\sum_{u=0}^{k}\bigl(2^{2u-1}-1\bigr)\biggl(\begin{matrix}2k+1\\2u\end{matrix}\biggr)B_{2u}=0, \label{eq:recurrent_bern_1}
\end{equation}
$$
\sum_{u=0}^{k}2^{2u-1}\biggl(\begin{matrix}2k+1\\2u\end{matrix}\biggr)B_{2u}=\sum_{u=0}^{k}\biggl(\begin{matrix}2k+1\\2u\end{matrix}\biggr)B_{2u}=\frac{2k+1}{2} \;\;\;\;\; \forall k,
$$
\begin{equation}
\sum_{u=0}^{k}\bigl(2^{2u-1}-1\bigr)\biggl(\begin{matrix}2k\\2u\end{matrix}\biggr)B_{2u}=-2^{2k-1}B_{2k} \;\;\;\;\; \forall k, \label{eq:recurrent_bern_2}
\end{equation}
$$
\sum_{u=0}^{(k-1)/2}\frac{\biggl(\begin{matrix}k\\2u+1\end{matrix}\biggr)}{k-u}\:B_{2(k-u)}=\int\limits_{0}^{-1}(x^2+x)^{k}dx \;\;\;\;\; \forall k,
$$
$$
\sum_{u=0}^{k}\:\frac{\biggl(\:\begin{matrix}k\\u\end{matrix}\:\biggr)}{2k+1-u}\biggl(B_{2k+1-u}+\frac{(-1)^u}{2}\biggr)=0 \;\;\;\;\; \forall k,
$$
$$
\sum_{u=0}^{2k-1}\:\frac{\biggl(\begin{matrix}2k+1\\u\end{matrix}\biggr)}{2^u}\biggl(B_{2k+1-u}+\frac{(-1)^u}{2}\biggr)=0 \;\;\;\;\; \forall k,
$$
or, in a more general form,
$$
\sum_{\begin{matrix}u_0+u_1+...+u_t=m\\u_i\geq 0\end{matrix}}\frac{k^{u_0}\prod\limits_{i=1}^{t-1}\biggl(\biggl(\:\begin{matrix}k\\2i\end{matrix}\:\biggr)B_{2i}\biggr)^{u_i}}{2^{u_0}\:\theta(u)\:u_0!\:u_1!...u_t!}\biggl(B_{\theta(u)}+\frac{(-1)^{u_0}}{2}\biggr)=0 \;\;\; \forall k, m \; (k>1),
$$
where $t=[(k+1)/2]$, $2\,|\,km$ and $\theta(u)=km+1-u_0-2\sum\limits_{i=1}^{t-1}tu_i$. 

And also we have
\begin{equation}
\sum_{u=0}^{2m+1}\frac{(-1)^{u-1}\:a_u}{n+1-u}=\sum_{u=0}^{(n-2)/2}\frac{a_{2u+1}}{m+1-u}\:B_{n-2u}, \label{eq:recurrent_bern_3}
\end{equation}
where $n$ is an integer, $a_0, a_2,..., a_{2m}$, $m=[(n-1)/2]$, is a sequence of real numbers, and
$$
a_{2k+1}=\sum_{\nu=0}^{k}\:\frac{(n-2\nu)!}{(n-2k-1)!}\:a_{2\nu}\sum_{\sigma=0}^{k-\nu}\:\frac{(-1)^{\sigma}}{2^{1+\sigma}\:u_0!(2u_1)!...(2u_{\sigma}!)},
$$
where $u_0+u_1+...+u_{\sigma}=2k+1-2\nu$, $u_i>0$ and $0\leq k\leq m$.
\begin{remark}
As a special case of formula (\ref{eq:recurrent_bern_3}), putting $a_{2u}=\biggl(\begin{matrix}n\\2u\end{matrix}\biggr)B_{2u}$, $u=0, 1,...,m$, we get
$$
\sum_{u=0}^{2m+1}\:\frac{(-1)^{u-1}\biggl(\:\begin{matrix}n\\u\end{matrix}\:\biggr)B_u}{n+1-u}=B_n,
$$
or, which is the same,
$$
\sum_{u=0}^{n}(-1)^{u-1}\biggl(\begin{matrix}n+1\\u\end{matrix}\biggr)B_u=0
$$
that is the Louivre's formula.
\end{remark}

\subsection{Finite Formulas of Bernoulli Numbers}

Analogously, the finite (non-recurrence) formulas of Bernoulli numbers can be derived
\begin{eqnarray}
	B_k = k!\sum_{t=1}^{k}\frac{(-1)^{t+k}}{(u_1+1)!(u_2+1)!...(u_t+1)!}\;\;\;\;\;\;\;\;\; \forall k,  \label{eq:finite_bern_1}   \\
	B_{2k}=\frac{(2k)!}{2}\sum_{t=1}^{k}\frac{(-1)^{t-1}(2u_1-1)}{(2u_1+1)!(2u_2+1)!...(2u_t+1)!} \;\;\;\;\; \forall k,  \label{eq:finite_bern_2}
\end{eqnarray}
where $u_1+u_2+...+u_t=k$, \: $u_t>0$,
\begin{eqnarray}
	B_{2k}=\frac{(2k)!}{2^{2k}}\sum_{t=0}^{k}\frac{(-1)^{t}}{(2u_0)!(2u_1+1)!(2u_2+1)!...(2u_t+1)!} \;\;\;\;\; \forall k, \label{eq:finite_bern_3}
\end{eqnarray}
where $u_0+u_1+u_2+...+u_t=k$, \: $u_0\geq 0$, $u_i>0$, $i>0$,
\begin{eqnarray}	B_{2k}=\frac{1}{2^{2k+1}-2}\biggl(1+(2k)!\sum_{u=1}^{k-1}\sum_{\sigma=1}^{k-u}\frac{(-1)^{\sigma}}{2^{\sigma}(2u-1)! \: \nu_1!\nu_2!...\nu_{\sigma}!}\biggr) \; \forall k  \label{eq:finite_bern_4}
\end{eqnarray}
where $\nu_1+\nu_2+...+\nu_{\sigma}=2k+1-2u$, \: $1<\nu_i\leq \nu_{i+1}$, and
\begin{eqnarray}
	B_k = (-1)^{k}\sum_{u=1}^{k}\frac{1}{u+1}\sum_{t=1}^{u}(-1)^{t}\:t^{k}\biggl(\:\begin{matrix}u\\ t\end{matrix}\:\biggr)  \forall k.  \label{eq:finite_bern_5}
\end{eqnarray}

\begin{example}
Above formulas can be easily checked. Putting $k=3$ in (\ref{eq:finite_bern_1}) and (\ref{eq:finite_bern_2}), we have
$$
B_3 = 6\sum_{t=1}^{3}\frac{(-1)^{t-1}}{(u_1+1)!...(u_t+1)!}
$$
and
$$
B_6 = 360\sum_{t=1}^{3}\frac{(-1)^{t-1}(2u_1-1)}{(2u_1+1)!...(2u_t+1)!},
$$
where $u_1+...+u_t=3$ \: $(u_i>0)$, that is $u_1=3$ for $t=1$, $u_1=1$, $u_2=2$ and $u_1=2$, $u_2=1$ for $t=2$ and $u_1=1$, $u_2=1$, $u_3=1$ for $t=3$. Therefore
$$
B_3 = 6 \biggl(\frac{1}{4!}-\frac{2}{2!\:3!}+\frac{1}{2!\:2!\:2!}\biggr)=6\biggl(\frac{1-4+3}{24}\biggr)=0,
$$
i. e.
$$
B_3=0
$$
and
$$
B_6 = 360\biggl(\frac{5}{7!}-\frac{4}{3!\:5!}+\frac{1}{3!\:3!\:3!}\biggr)=360\biggl(\frac{15-84+70}{24}\biggr)=\frac{1}{42},
$$
i. e.
$$
B_6=\frac{1}{42}.
$$
Putting $k=2$ in (\ref{eq:finite_bern_3}) and (\ref{eq:finite_bern_4}), we get
$$
B_4 = \frac{3}{2}\sum_{t=0}^{2}\frac{(-1)^{t}}{(2u_0)!(2u_1+1)!...(2u_t+1)!},
$$
where $u_0+u_1+...+u_t=2, \: (u_0\geq 0, u_i>0, t\geq 1)$, that is $u_0=2$ for $t=0$, $u_0=0$, $u_1=2$ and $u_0=1$, $u_1=1$ for $t=1$ and $u_0=1$, $u_1=1$, $u_2=1$ for $t=2$. Then
$$
B_4 = \frac{3}{2}\biggl(\frac{1}{4!}-\biggl(\frac{1}{5!}+\frac{1}{2!\:3!}\biggr)+\frac{1}{3!\:3!}\biggr)= \frac{3}{2}\biggl(\frac{15-3-30+10}{360}\biggr)=-\frac{1}{30},
$$
$$
B_4=-\frac{1}{30}
$$
and
$$
B_4 = \frac{1}{30}\biggl(1+24\sum_{t=1}^{1}\frac{-1}{2u_1!}\biggr),
$$
where $u_1=3$, i. e.
$$
B_4 = \frac{1}{30}\biggl(1+24\biggl(\frac{-1}{12}\biggr)\biggr)=-\frac{1}{30},
$$
$$
B_4 = -\frac{1}{30}.
$$
And according to (\ref{eq:finite_bern_5}) for $k=2, 3$, we have
$$
B_2 = \frac{1}{2}(-1)+\frac{1}{3}(-2+4)=-\frac{1}{2}+\frac{2}{3}=\frac{1}{6},
$$
$$
B_2=\frac{1}{6}
$$
and
$$
B_3 = \frac{1}{2}(-1)+\frac{1}{3}(-2+8)+\frac{1}{4}(-3+8\cdot3-27)=-\frac{1}{2}+2-\frac{3}{2}=0,
$$
$$
B_3=0.
$$
\end{example}

Using formulas (\ref{eq:rem_qu_zero}), (\ref{eq:2_theor_gen_inf_sum}), (\ref{sigma_fraction}), (\ref{eq:analysis1}) and (\ref{eq:analysis2}), we have
\begin{eqnarray}
\int\limits_{-1}^{0}Q_r(x)dx=0, \;\; N_r\equiv1\;(\!\!\!\!\!\!\mod2), \label{eq:qu_1}   \\
\int\limits_{-1}^{0}Q_r(x)dx=-2\sum_{u=1}^{N_r}\frac{q_u B_{u+1}}{u+1}, \;\; N_r\equiv0\;(\!\!\!\!\!\!\mod2), \label{eq:qu_2}   \\	\sum_{t=1}^{k-1}\:\sum_{\begin{matrix}u_0+u_1+u_t=k\\(u_i>0)\end{matrix}}\:\frac{(-1)^{t}\:2^{1-t}(B_{2u_1}-1/2)}{(2u_0-1)!(2u_1)!(2u_2)!...(2u_t)!}=\frac{1}{(2k-1)!} \;\;\;\;\;\;\;\;\; \forall k>1, \label{eq:qu_3} \\
\sum_{u=1}^{\infty}(2u-1)^{k-1}=\frac{2\:(2^{k-1}-1)B_k+(-1)^k}{2k} \;\;\;\;\;\;\; \forall k,  \label{eq:qu_4}  \\
\sum_{u=1}^{\infty}(-1)^{u-1}(2u-1)^{k-1}=\frac{-1}{2k}\sum_{u=1}^{k}(-1)^{u}\:2^{u}(2^{u}-1)\biggl(\:\begin{matrix}k\\ u\end{matrix}\:\biggr)B_{u} \; \forall k \label{eq:qu_5}
\end{eqnarray}
or, in view of (\ref{eq:from_diver_sin})
$$
\sum_{u=1}^{\infty}(4u-3)^{k-1}=\frac{1}{4k}\bigl((-3)^{k}+(2^{k}-2)B_k\bigr)-\sum_{u=1}^{k}(-1)^{u}\:2^{u}(2^{u}-1)\biggl(\:\begin{matrix}k\\ u\end{matrix}\:\biggr)B_{u}
$$
and
$$
\sum_{u=1}^{\infty}(4u-3)^{2k-1}=\frac{1}{8k}\bigl(3^{2k}+(2^{2k}-2)B_{2k}\bigr).
$$

\begin{example}
Let $Q_2(x)=18B_2(x)B_5(x)$. Then
$$
Q_2(x)=\biggl(x^3+\frac{3}{2}x^2+\frac{1}{2}x\biggr)\biggl(x^6+3x^5+\frac{5}{2}x^4-\frac{1}{2}x^2\biggr)
$$
i. e.
$$
Q_2(x)=x^9+\frac{9}{2}x^8+\frac{15}{2}x^7+\frac{21}{4}x^6+\frac{3}{4}x^5-\frac{3}{4}x^4-\frac{1}{4}x^3, \;\;\;\;\; \bigl(9\equiv1\;(\!\!\!\!\!\!\mod2)\bigr)
$$
and
$$
\int\limits_{-1}^{0}Q_2(x)dx=-\frac{1}{10}+\frac{1}{2}-\frac{15}{16}+\frac{3}{4}-\frac{1}{8}-\frac{3}{20}+\frac{1}{16}=0,
$$
which confirms the formula (\ref{eq:qu_1}).

Now let $Q_2(x)=36\bigl(B_2(x)\bigr)^{2}B_3(x)$. Then
$$
Q_2(x)=\biggl(x^3+\frac{3}{2}x^2+\frac{1}{2}x\biggr)^{2}\bigl(x^4+2x^3+x^2\bigr)
$$
i. e.
$$
Q_2(x)=x^{10}+5x^9+\frac{41}{4}x^8+11x^7+\frac{13}{2}x^6+2x^5+\frac{1}{4}x^4, \;\;\;\;\; \bigl(10\equiv0\;(\!\!\!\!\!\!\mod2)\bigr).
$$
Then, on the one hand
$$
\int\limits_{-1}^{0}Q_2(x)dx=\frac{1}{11}-\frac{1}{2}+\frac{41}{36}-\frac{11}{8}+\frac{13}{14}-\frac{1}{3}+\frac{1}{20}=\frac{1}{27720},
$$
from the other hand
$$
-\sum_{u=1}^{10}\frac{q_u\:B_{u+1}}{u+1}=-\frac{1}{2}B_{10}-\frac{11}{8}B_{8}-\frac{1}{3}B_{6}=-\frac{5}{132}+\frac{11}{240}-\frac{1}{126}=\frac{1}{55440}=\frac{1}{2}\Bigl(\frac{1}{27720}\Bigr).
$$
Comparing these two, we get
$$
\int\limits_{-1}^{0}Q_2(x)dx=-2\sum_{u=1}^{10}\frac{q_u\:B_{u+1}}{u+1},
$$
that is the formula (\ref{eq:qu_2}).

And now let us consider the equality (\ref{eq:qu_3}) for $k=4$. Then we have
\begin{equation}
L = \sum_{t=1}^{3}\frac{(-1)^{t}\:2^{1-t}\bigl(B_{2u_1}-1/2\bigr)}{(2u_0-1)!(2u_1)!...(2u_t)!}, \label{eq:qu_qu_3}
\end{equation}
where $u_0+u_1+...+u_t=4$, \; $u_i>0$, i. e.
\begin{tabbing}

22222 \= 2222222222222222 \= 2222222222222222222 \= 222222222222  \kill

\>  for $t=1$        \> for $t=2$             \> for $t=3$          \\
\> \underline{$u_0+u_1=4$} \> \underline{$u_0+u_1+u_2=4$} \> \underline{$u_0+u_1+u_2+u_3=4$}  \\ 
\> $1+3=4$        \> $1+1+2=4$             \> $1+1+1+1=4$           \\
\> $2+2=4$          \> $1+2+1=4$            \>                 \\ 
\> $3+1=4$         \> $2+1+1=4$             \>                  

\end{tabbing} 
Substituting these values of $u_i$ into (\ref{eq:qu_qu_3}), we get
\begin{eqnarray}
	L &=& -\biggl(\frac{B_6-1/2}{6!}+\frac{B_4-1/2}{3!\:4!}+\frac{B_2-1/2}{5!\:2!}\biggr)+ 
	\frac{1}{2}\biggl( \frac{B_2-1/2}{2!\:4!}+\frac{B_4-1/2}{4!\:2!}+ \nonumber \\
	&&+
	\frac{B_2-1/2}{3!\:2!\:2!}\biggr)-\frac{1}{4}\biggl(\frac{B_2-1/2}{2!\:2!\:2!}\biggr)=
	\frac{10}{6!\:21}+\frac{8}{3!\:4!\:15}+\frac{1}{5!\:2!\:3}- \nonumber \\
	&&-
	\frac{1}{2}\biggl(\frac{1}{2!\:4!\:3} - \frac{8}{4!\:2!\:15}-\frac{1}{3!\:2!\:2!\:3}\biggr)+
	\frac{1}{4}\biggl(\frac{1}{2!\:2!\:2!\:3}\biggr)= \frac{1}{7!}. \nonumber 
\end{eqnarray}
Therefore
$$
\sum_{t=1}^{3}\frac{(-1)^{t}\:2^{1-t}\bigl(B_{2u_1}-1/2\bigr)}{(2u_0-1)!(2u_1)!...(2u_t)!}= \frac{1}{7!},
$$
where $u_0+u_1+u_t=4$, which is in agreement with (\ref{eq:qu_3}).
\end{example}

\subsection{Theorems Involving Bernoulli Numbers}

Let $f(x)=\sum_{u=1}^{2k}b_u x^u$ be a polynomial of even degree satisfying the condition (\ref{eq:theor_vartheta}).
\begin{theorem} \label{theor:bern_sum}
With the above notation of Theorem \ref{theor:gen_inf_sum},
\begin{equation} \label{eq:theor_bern_sum}
\sum_{u=1}^{k}\frac{b_{2u-1}}{u}B_{2u}=-\epsilon \: \int\limits_{-1}^{0}\biggl(\sum_{u=\delta}^{t-1+\delta}\bigl(f(x-\epsilon u)-f(-\epsilon u)\bigr)\biggr)dx.
\end{equation}
\end{theorem}
\begin{proof}
According to (\ref{eq:gen_inf_sum}) and (\ref{eq:analysis1}), we have
$$
\sum_{u=1}^{\infty}f(u)=-\sum_{u=1}^{2k}\frac{b_u}{u+1}B_{u+1}=\frac{\epsilon}{2}\sum_{u=\delta}^{t-1+\delta}\bigl(\lim_{n\rightarrow\infty}f(n-\epsilon u)-f(-\epsilon u)\bigr)-\frac{1}{2}f(0).
$$
But since $f(0)=0$ and $B_{2u-1}=0$, $u>1$, in view of (\ref{eq:theor_polynom2})
$$
\sum_{u=1}^{k}\frac{b_{2u-1}}{u}B_{2u}=-\epsilon \: \int\limits_{-1}^{0}\biggl(\sum_{u=\delta}^{t-1+\delta}\bigl(f(x-\epsilon u)-f(-\epsilon u)\bigr)\biggr)dx.
$$
The theorem is proved.
\end{proof}

\begin{example}
The polynomial $f(x)=x^4+6x^3+x^2-24x$ satisfies the condition (\ref{eq:theor_vartheta}) for $t=3$ and $\epsilon=1$. Indeed,
$$
f(-x)=x^4-6x^3+x^2+24x,
$$
\begin{eqnarray}
	f(x-3)\!\!\!&=&\!\!\!(x^4-12x^3+54x^2-108x+81)+(6x^3-54x^2+162x-162)+ \nonumber \\
	&&+
	(x^2-6x+9)-(24x-72)=x^4-6x^3+x^2+24x \nonumber
\end{eqnarray}
and therefore
$$
f(-x)=f(x-3).
$$
Then, we get
\begin{eqnarray}	
&&\sum_{u=0}^{2}\bigl(f(x-u)-f(-u)\bigr)\!\!\!= \!\!\!(x^4+6x^3+x^2-24x)+\bigl((x-1^{4})+6(x-1)^{3}+ \nonumber \\
&&+
(x-1)^{2}-24(x-1)\bigr)+ \bigl((x-2)^{4}+6(x-2)^{3}+  
(x-2)^{2}-24(x-2)\bigr)= \nonumber \\
&&=
3x^4+6x^3-21x^2-24x \nonumber
\end{eqnarray}
and
\begin{eqnarray}	\int\limits_{-1}^{0}\sum_{u=0}^{2}\bigl(f(x-u)-f(-u)\bigr)dx&=&\int\limits_{-1}^{0}(3x^4+6x^3-21x^2-24x)dx= 
\frac{41}{10}. \nonumber
\end{eqnarray}
And since
$$
\sum_{u=1}^{2}\frac{b_{2u-1}}{u}B_{2u}=-\frac{24}{1}B_{2}+\frac{6}{2}B_{4}=-4-\frac{1}{10}=
-\frac{41}{10},
$$
hence
$$
\sum_{u=1}^{2}\frac{b_{2u-1}}{u}B_{2u}=\int\limits_{-1}^{0}\sum_{u=0}^{2}\bigl(f(x-u)-f(-u)\bigr)dx,
$$
which Theorem \ref{theor:bern_sum} do asserts.
\end{example}

\begin{remark}
If a polynomial of even degree $\nu$ satisfies the condition
\begin{equation} \label{eq:polycond}
f(-x)=f(x-1)\: \footnote{\:It is clear that the coefficient at $x^{\nu-1}$ of $f(x)$ is equal to $\nu/2$} 
\end{equation}
then the polynomial $f^{(t)}(x)=t^{\nu}f(x/t)$ satisfies the equality
$$
f^{(t)}(-x)=f^{(t)}(x-t).
$$
Hence, taking a monic polynomial $f(x)$ of degree $\nu$ with integer coefficients satisfying the condition (\ref{eq:polycond}), we have
$$
f^{(t)}(n)=n^{\nu}+\frac{t\nu}{2}n^{\nu-1}+t^{2}N(n),
$$
where $N(n)$ is an integer.

And since for any natural number $q$
$$
q\sum_{u=0}^{t-1} \: \int\limits_{-1}^{0}\bigl((x-u)^{q-1}-(-u)^{q-1})\bigr)dx=(-1)^{q-1} \: t\sum_{u=0}^{q-2}S_{u}^{t}\biggl(\:\begin{matrix}q\\ u\end{matrix}\:\biggr),
$$
where $S_{u}^{t}=\frac{1}{t}\sum_{i=0}^{t-1}i^{u}$, then putting $t=p$ and $\theta=q$, in view of (\ref{eq:theor_bern_sum}), we obtain the following
\begin{corollary}[Theorem \ref{theor:bern_sum}]
For any prime number $p>2$ and odd \: $\theta$, $1<\theta\leq p$, the following congruence holds
\begin{equation}
\theta B_{\theta-1}\equiv -\sum_{u=0}^{\theta-2}S_{u}^{p}\biggl(\:\begin{matrix}\theta\\ u\end{matrix}\:\biggr) \;\; (\!\!\!\!\!\!\mod p). \label{eq:bern_prime_cong}
\end{equation}
Whence, in particular, $pB_{p-1}\equiv -1 \;\; (\!\!\!\!\mod p)$.
\end{corollary}
\end{remark}


Let $\theta_{p}(x)=\sum_{u=1}^{p-1}a_{u}x^{u}$ be a polynomial of degree $p-1$, $p\geq 1$, with integer coefficients, where $a_{2u-1}\equiv 0\;(\!\!\!\!\mod p)$, $u=\overline{1, \, 2^{-1}(p-9)}$. And let $\theta_{p}(x)$ satisfies the condition (\ref{eq:polycond}). Then

\begin{proposition} \label{prop:fermat}
If \: $a_{2}\equiv \!\!\!\!\!\!/ \;\;0\;(\!\!\!\!\mod p)$, then the first case of Fermat's theorem holds true for prime power $p\geq 11$.
\end{proposition}
\begin{proof}
Indeed,  the polynomial $H(x)=\theta_{p}(x)B_{p-3}(x)$ in view of (\ref{eq:1_from_lemma3}) satisfies the condition (\ref{eq:1_theor_theta1}) for $\epsilon=1$. This means that, according to (\ref{eq:2_theor_theta1}, the following congruence holds
$$
p\:\int\limits_{0}^{-1}H(x)dx \equiv \frac{a_{2}}{2}+\sum_{u=1}^{4}u\,a_{p-2u}\,B_{p-2u-1} 
\equiv 0 \:\:(\!\!\!\!\!\!\mod p).
$$
Whence, from the Mirimanov's theorem, which asserts that if the numerator of at least one of the four Bernoulli numbers $B_{p-3}, B_{p-5}, B_{p-7}, B_{p-9}$ is not divided by $p$, then the first case of the Fermat's theorem holds true for prime exponent $p$, the proposition follows.
\end{proof}

\begin{remark}
By not complicated calculations it is established that for any prime number $p\geq 11$ there exists a polynomial $\theta_{p}(x)$ satisfying the condition (\ref{eq:polycond}), which coefficient $a_2$ at $x^2$ satisfies the congruence
$$
a_2 \equiv \sum_{t=1}^{q}(-1)^{t}\:T_{q+2}^{2u_{0}+5}\:T_{q-u_0}^{2u_1}\:...\:T_{q-u_0-...-u_{t-1}}^{2u_t}\;(\bmod\;p),
$$
where $q=\frac{p-5}{2}$, $u_0+u_1+...+u_t=q+2$, $u_i>0$, and $T_{n}^{m}==\frac{(n-1)!}{n!(n-m)!}$.
Hence, the first case of Fermat's theorem holds true for any prime exponent $p\geq 11$ if
$$
\sum_{t=1}^{q}(-1)^{t}\:T_{q+2}^{2u_{0}+5}\:T_{q-u_0}^{2u_1}\:...\:T_{q-u_0-...-u_{t-1}}^{2u_t}\equiv\!\!\!\!\! /\;\;0\;(\bmod\;p).
$$

\subsection{Some Results Involving Symmetric Functions}

It is easy to verify that one of the polynomials of degree $2k$ satisfying the condition (\ref{eq:polycond}) is the following one
$$
f(x)=\sum_{u=0}^{t}\sigma_{u}(t)\:x^{t+1-u},
$$
where $\sigma_{0}(t)=1$, $\sigma_{u}(t)=(1\cdot2\cdot\cdot\cdot u)+...+\bigl((t+1-u)\cdot\cdot\cdot t\bigr)$, $u>0$, is an elementary symmetric function of natural numbers of the interval $|1, t|$, $t=2k-1$. In a more general form,
$$
f(x)=\sum_{u=0}^{k}\sum_{\nu=0}^{k}(-1)^{u}\:a_{u}\:\overline{a_{\nu}}\:x^{u+\nu},
$$
where $a_{u}$ are any real numbers, and the coefficients $\overline{a_{\nu}}$ are defined by the formula
\begin{equation}
\overline{a_{\nu}}=\sum_{u=\nu}^{k}(\epsilon t)^{u}\:a_{u}\biggl(\:\begin{matrix}u\\ \nu\end{matrix}\:\biggr). \label{eq:symm_1}
\end{equation}
Then, according to the Theorem \ref{theor:bern_sum}, we get
\begin{equation}
\sum_{u=0}^{t}\frac{\sigma_{u}(t)}{t+2-u}\:B_{t+2-u}=-\frac{1}{2}\int\limits_{-1}^{0}\biggl(\sum_{u=0}^{t-1}f(x-u)\biggr)dx  \label{eq:symm_2}
\end{equation}
and
\begin{equation}
\sum_{u=0}^{k}\sum_{\nu=0}^{k}(-1)^{u}\:\frac{a_{u}\:\overline{a_{\nu}}}{u+\nu+1}\biggl(B_{u+\nu+1}+\frac{(-1)^{u+\nu}}{2}\biggr)=a_{0}\:\overline{a_{0}}, 
\;\;\;\;\; (\epsilon t=1). \label{eq:symm_3}
\end{equation}
The analogous equality immediately follows from the formula (\ref{eq:theor_binom}, namely
\begin{equation}
\sum_{u=1}^{(k+1)/2}\sigma_{k+1-2u}(k)\:\frac{2^{2u}-1}{u}\:B_{2u}=k!\:(1-\frac{1}{2^k}) \;\;\;\;\;\;\; \forall k. \label{eq:symm_4}
\end{equation}
Moreover, for any natural number $m$ and real $\alpha$ the following two formulas hold
$$
\sum_{\begin{matrix}u_0+u_1+...+u_{\sigma}=m\\ u_i\geq 0\end{matrix}}\;\frac{B_{\tau(u)}\prod\limits_{i=1}^{k-1}\biggl(\:\begin{matrix}k\\ i\end{matrix}\:\biggr)^{u_i+u_{k+i}}}{\tau(u)\:\theta(u)\:u_{0}!\:u_{1}!\:...\:u_{\sigma}!}\:\biggl(B_{\theta(u)}-\frac{(-1)^{\theta(u)}}{2}\biggr)=0, \;\;\;\;\; 2\:|\!\!/\:k
$$
and
$$
\sum_{\begin{matrix}u_0+u_1+...+u_{\sigma}=m\\ u_i\geq 0\end{matrix}}\frac{\alpha^{\tau(u)}(2\alpha+1)^{u_{k}}\prod\limits_{i=1}^{k-1}\biggl(\:\begin{matrix}k\\i\end{matrix}\:\biggr)^{u_i+u_{k+i}}}{\theta(u)\:u_{0}!\:u_{1}!\:...\:u_{\sigma}!}\:\biggl(B_{\theta(u)}-\frac{(-1)^{\theta(u)}}{2}\biggr)=0, \;\;\; 2\:|\:k,
$$
where $\sigma=2k-1$, \: $\tau(u)=1+\sum\limits_{i=1}^{k-1}u_{k+i}$ \: and \: $\theta(u)=2km+1-\sum\limits_{i=1}^{\sigma}i\:u_{i}$, assuming $\prod\limits_{u=1}^{0}\cdot=1$ and $\sum\limits_{u=1}^{0}\cdot=0$.
\end{remark}

\begin{example}
Let us check two above formulas (\ref{eq:symm_2}) and (\ref{eq:symm_3}). Putting $t=3$ in (\ref{eq:symm_2}), we get
$$
\sum_{u=0}^{3}\frac{\sigma_{u}(3)}{5-u}\:B_{5-u}=\frac{6}{4}B_4+\frac{6}{2}B_2=-\frac{1}{20}+
\frac{1}{2}=\frac{9}{20}
$$
and since $f(x)=x^4+6x^3+11x^2+6x$ and $\sum\limits_{u=0}^{2}f(x-u)=3x^4+6x^3+9x^2+6x$, we have
$$
\int\limits_{-1}^{0}\biggl(\sum_{u=0}^{2}f(x-u)\biggr)dx=\int\limits_{-1}^{0}\bigl(3x^4+6x^3+9x^2+6x\bigr)dx=\frac{3}{5}-\frac{3}{2}+3-3=-\frac{9}{10}.
$$
Hence
$$
\sum_{u=0}^{3}\frac{\sigma_{u}(3)}{5-u}\:B_{5-u}=-\frac{1}{2}\int\limits_{-1}^{0}\biggl(\sum_{u=0}^{2}f(x-u)\biggr)dx,
$$
which confirms the formula (\ref{eq:symm_2}).

Now putting in (\ref{eq:symm_3}) $k=1$, $a_0=\frac{1}{2}$, $a_1=\frac{1}{3}$, and according to (\ref{eq:symm_1}) \: $\overline{a_0}=\frac{5}{6}$ and $\overline{a_1}=\frac{1}{3}$, we get
\begin{eqnarray}	&&\sum_{u=0}^{1}\sum_{\nu=0}^{1}\!\!(-1)^{u}\!\!\frac{a_{u}\:\overline{a_{\nu}}}{u+\nu+1}\biggl(B_{u+\nu+1}-\frac{(-1)^{u+\nu}}{2}\biggr)\:= 
\:\frac{a_{0}\:\overline{a_0}}{1}\Bigl(B_1+\frac{1}{2}\Bigr)+ \nonumber \\
&&+
\!\frac{a_{0}\:\overline{a_1}}{2}\Bigl(B_1-\frac{1}{2}\Bigr)- 
\frac{a_{1}\:\overline{a_0}}{2}\Bigl(B_2-\frac{1}{2}\Bigr)- 
\frac{a_{1}\:\overline{a_1}}{3}\Bigl(B_3+\frac{1}{2}\Bigr)= 
\!\frac{5}{12}\Bigl(B_1+\frac{1}{2}\Bigr)+ \nonumber \\
&&+
\frac{1}{12}\Bigl(B_2-\frac{1}{2}\Bigr)-
\frac{5}{36}\Bigl(B_2-\frac{1}{2}\Bigr)-\frac{1}{27}\Bigl(B_3+\frac{1}{2}\Bigr)= 
\frac{5}{12}=a_0\:\overline{a_0}, \nonumber
\end{eqnarray}
which is in agreement with (\ref{eq:symm_3}).
\end{example}

\vspace{1cm}
\begin{center}
\rule{9.5cm}{0.25mm}
\end{center}

\newpage

\section{Trigonometric Functions and Series}\label{sec:trig}

All the trigonometric functions are elementary functions. So, according to the Remark 2.6., these functions have a certain limit on the set $\mathbb{Z}$.

Below we write out several algebraic relations for some trigonometic functions.
\begin{enumerate}
\item $\sin\alpha \pm \sin\beta=2\sin\Bigl(\frac{\alpha\pm\beta}{2}\Bigr)\cos\Bigl(\frac{\alpha\mp\beta}{2}\Bigr)$  \label{item:1}  \\
\item $\sin(\alpha \pm \beta)=\sin\alpha\cos\beta \pm \cos\alpha\sin\beta$     \label{item:2}   \\
\item $\cos\alpha-\cos\beta=
-2\sin\Bigl(\frac{\alpha-\beta}{2}\Bigr)\sin\Bigl(\frac{\alpha+\beta}{2}\Bigr)$  \label{item:3}  \\
\item $\cos\alpha+\cos\beta=
2\cos\Bigl(\frac{\alpha-\beta}{2}\Bigr)\cos\Bigl(\frac{\alpha+\beta}{2}\Bigr)$  \label{item:4} \\
\item $\cos(\alpha \pm \beta)=\cos\alpha\cos\beta \mp \sin\alpha\sin\beta$ \label{item:5}  \\
\item $\tan\alpha \pm \tan\beta=\frac{\sin(\alpha \pm \beta)}{\cos\alpha\cos\beta}$ \label{item:6} \\
\end{enumerate}

According to the item (\ref{item:6}),
$$
\tan(n+\frac{1}{2})\theta-\tan(n-\frac{1}{2})\theta=\frac{\sin\theta}{\cos(n+\frac{1}{2})\theta \, \cos(n-\frac{1}{2})\theta},
$$
where $-\pi<\theta<\pi$, and in view of axiom (\ref{axiom:5})
$$
\tan(n+\frac{1}{2})\theta-\tan\frac{1}{2}\theta=\sum_{u=1}^{n}\frac{\sin\theta}{\cos(u+\frac{1}{2})\theta \, \cos(u-\frac{1}{2})\theta} \;\;\;\;\; \forall n.
$$
Now passing to the limit and taking into account axiom (\ref{axiom:1}) and formula (\ref{eq:inf_sum}), we get
$$
\lim_{n\rightarrow\infty}\tan(n+\frac{1}{2})\theta-\tan\frac{1}{2}\theta=\sum_{u=1}^{\infty}\frac{\sin\theta}{\cos(u+\frac{1}{2})\theta \, \cos(u-\frac{1}{2})\theta}=-\tan\frac{1}{2}\theta
$$
and
$$
\lim_{n\rightarrow\infty}\tan(n+\frac{1}{2})\theta=0 \;\;\;\;\;\;\; -\pi<\theta<\pi.
$$
Analogously we also obtain the equality
$$
\lim_{n\rightarrow\infty}\cot(n+\frac{1}{2})\theta=0 \;\;\;\;\;\;\; 0<\theta<2\pi.
$$

\begin{remark}
The last two equalities are, in fact, the special cases of the formula (\ref{eq:3_cor_lemma2}).
\end{remark}

According to item (\ref{item:1}), for any real number $\theta$
$$
\sin(2n+1)\theta - \sin\theta = 2\sin n\theta\cos(n+1)\theta
$$
and
$$
\sin(2n+1)\theta + \sin\theta = 2\sin(n+1)\theta\cos n\theta.
$$
Passing to the limit and taking into account the formula (\ref{eq:3_cor_lemma2}), we get
$$
\lim_{n\rightarrow\infty}\sin n\theta \cos(n+1)\theta=-\frac{\sin\theta}{2}
$$
and
$$
\lim_{n\rightarrow\infty}\sin(n+1)\theta \cos n\theta=\frac{\sin\theta}{2}.
$$
Relying on item (\ref{item:2}), we have
$$
\sin\bigl((n+\frac{1}{2}) \pm \frac{1}{2}\bigr)\theta=\sin(n+\frac{1}{2})\theta\cos\frac{1}{2}\theta \pm \cos(n+\frac{1}{2})\theta\sin\frac{1}{2}\theta
$$
and passing to the limit with use of (\ref{eq:3_cor_lemma2}), we get
$$
\lim_{n\rightarrow\infty}\sin(n+1)\theta=\sin\frac{1}{2}\theta\lim_{n\rightarrow\infty}\cos(n+\frac{1}{2})\theta,
$$
$$
\lim_{n\rightarrow\infty}\sin n\theta = -\sin\frac{1}{2}\theta\lim_{n\rightarrow\infty}\cos(n+\frac{1}{2})\theta
$$
and 
$$
\lim_{n\rightarrow\infty}\sin(n+1)\theta=-\lim_{n\rightarrow\infty}\sin n\theta.
$$
The last equality is in good agreement with Lemma \ref{lemma2}.

In analogous way we also derive the equality
$$
\lim_{n\rightarrow\infty}\cos(n+1)\theta=\lim_{n\rightarrow\infty}\cos n\theta,
$$
which as well is in good agreement with Lemma \ref{lemma2}.

According to item (\ref{item:1})
$$
\sin(n+\frac{1}{2})\theta-\sin(n-\frac{1}{2})\theta=2\sin\frac{1}{2}\theta\cos n\theta
$$
whence, passing to the limit, we get
\begin{equation}
\lim_{n\rightarrow\infty}\sin(n-\frac{1}{2})\theta=-2\sin\frac{1}{2}\theta\lim_{n\rightarrow\infty}\cos n\theta. \label{eq:trig_1}
\end{equation}

Using the items (\ref{item:2}) and (\ref{item:5}), we have
$$
\sin\bigl((n+\frac{1}{2})+(\alpha-\frac{1}{2})\bigr)\theta=\sin(n+\frac{1}{2})\theta\cos(\alpha-\frac{1}{2})\theta+\cos(n+\frac{1}{2})\theta\sin(\alpha-\frac{1}{2})\theta
$$
and
$$
\cos\bigl((n+\frac{1}{2})+(\alpha-\frac{1}{2})\bigr)\theta=\cos(n+\frac{1}{2})\theta\cos(\alpha-\frac{1}{2})\theta-\sin(n+\frac{1}{2})\theta\sin(\alpha-\frac{1}{2})\theta,
$$
where $\alpha$ is a real number.

Now passing to the limit and using the formula (\ref{eq:3_cor_lemma2}), we obtain two equalities
$$
\lim_{n\rightarrow\infty}\sin(n+\alpha)\theta=\sin(\alpha-\frac{1}{2})\theta\lim_{n\rightarrow\infty}\cos(n+\frac{1}{2})\theta
$$
and
$$
\lim_{n\rightarrow\infty}\cos(n+\alpha)\theta=\cos(\alpha-\frac{1}{2})\theta\lim_{n\rightarrow\infty}\cos(n+\frac{1}{2})\theta,
$$
from which directly follows
$$
\cos(\alpha-\frac{1}{2})\theta\lim_{n\rightarrow\infty}\sin(n+\alpha)\theta= \sin(\alpha-\frac{1}{2})\theta\lim_{n\rightarrow\infty}\cos(n+\alpha)\theta.
$$
Particularly, for $\alpha=0$ and $-\pi<\theta<\pi$
\begin{equation}
\lim_{n\rightarrow\infty}\sin n\theta=-\tan\frac{1}{2}\theta\lim_{n\rightarrow\infty}\cos n\theta. \label{eq:trig_2}
\end{equation}

Let us take the odd function $\sin n\theta \cos n\theta$. Then, using the items (\ref{item:2}) and (\ref{item:5}), we have
\begin{eqnarray}
	&&\sin(n+\frac{1}{2})\theta \cos(n+\frac{1}{2})\theta=(\sin n\theta \cos\frac{1}{2}\theta + 
	\cos n\theta \sin\frac{1}{2}\theta)(\cos n\theta \cos\frac{1}{2}\theta - \nonumber \\
	&&-
	\sin n\theta \sin\frac{1}{2}\theta) = 
	\frac{1}{2}\cos^{2}\frac{1}{2}\theta\sin 2n\theta + \frac{1}{2}\sin\theta\cos^{2}n\theta - \frac{1}{2}\sin\theta \sin^{2}n\theta - \nonumber \\
	&&-
  \frac{1}{2}\sin^{2}\frac{1}{2}\theta \sin 2n\theta =  
  \frac{1}{2}\bigl((\cos^{2}\frac{1}{2}\theta - \sin^{2}\frac{1}{2}\theta)\sin 2n\theta + \sin\theta(\cos^{2}n\theta - 
  \sin^{2}n\theta)\bigr) =  \nonumber \\
  &&=
  \frac{1}{2}\bigl(\cos\theta\sin 2n\theta + \sin\theta(1-2\sin^{2}n\theta)\bigr)      \nonumber 
\end{eqnarray}
and, passing to the limit, we get
$$
0=-\sin\theta\lim_{n\rightarrow\infty}\sin^{2}n\theta + \frac{1}{2}\cos\theta \lim_{n\rightarrow\infty}\sin 2n\theta + \frac{\sin\theta}{2},
$$
i. e.
$$
\sin\theta\lim_{n\rightarrow\infty}\sin^{2}n\theta = \frac{1}{2}\cos\theta \lim_{n\rightarrow\infty}\sin 2n\theta + \frac{\sin\theta}{2}
$$
or, putting $0<|\theta|<\pi$ \; (in order $\sin\theta \neq 0$)

\begin{equation}
\lim_{n\rightarrow\infty}\sin^{2}n\theta = \frac{1}{2}\cot\theta \lim_{n\rightarrow\infty}\sin 2n\theta +\frac{1}{2}. \label{eq:trig_3}
\end{equation}

Relying on items (\ref{item:1}) and (\ref{item:3}), we have
\begin{eqnarray}
	&&\sin n\theta - \sin(n-1)\theta = 2\sin\frac{1}{2}\theta \cos(n-\frac{1}{2})\theta = 2\sin\frac{1}{2}\theta(\cos n\theta \cos\frac{1}{2}\theta + \nonumber \\
	&&+
	\sin n\theta \sin\frac{1}{2}\theta) = 
	\sin\theta \cos n\theta + 2\sin^{2}\frac{1}{2}\theta \sin n\theta \nonumber
\end{eqnarray}
i. e.
$$
\sin n\theta - \sin(n-1)\theta = \sin\theta \cos n\theta + 2\sin^{2}\frac{1}{2}\theta \sin n\theta
$$
and according to axiom (\ref{axiom:5})
$$
\sin n\theta = \sin\theta \sum_{u=1}^{n}\cos u\theta + 2 \sin^{2}\frac{1}{2}\theta \sum_{u=1}^{n} \sin u\theta \;\;\;\;\;\;\; \forall n.
$$
Now passing to the limit and taking into account (\ref{ex:cos}) and putting $0<\theta<2\pi$ \: (in order $\sin\frac{1}{2}\theta \neq 0$), we obtain
\begin{equation}
\sum_{u=1}^{\infty}\sin u\theta = \frac{1}{2}\cot\frac{1}{2}\theta + \frac{1}{2}\csc^{2}\frac{1}{2}\theta \lim_{n\rightarrow\infty}\sin n\theta. \label{eq:trig_4}
\end{equation}

According to item (\ref{item:3})
$$
\cos 2n\theta - \cos(2n-2)\theta = -2\sin\theta \sin(2n-1)\theta
$$
and in view of axiom (\ref{axiom:5}), we get the equality
$$
\cos 2n\theta - 1 = -2\sin^{2}\theta = -2\sin\theta \sum_{u=1}^{n}\sin(2u+1)\theta \;\;\;\;\;\;\; \forall n,
$$
which, when passing to the limit, gives
\begin{equation}
\sum_{u=1}^{\infty}\sin(2u-1)\theta = \csc\theta \lim_{n\rightarrow\infty}\sin^{2}n\theta, \;\;\;\;\;\;\;\;\; 0<|\theta|<\pi. \label{eq:trig_5}
\end{equation}


Now let us consider the expansion of $\sin x$ into power series
$$
\sin x = x-\frac{x^3}{3!}+\frac{x^5}{5!}-... \;\;\;\;\;\;\;\;\; |x|<\infty.
$$
Replacing $x$ by $x\theta$, we have
$$
\sin x\theta = x\theta-\frac{x^3 \theta^3}{3!}+\frac{x^5 \theta^5}{5!}-...,
$$
where $-\pi<\theta<\pi$. 

Applying the rule \ref{rule3} 
from Section \ref{sec:intro}, we get
$$
\theta \lim_{n\rightarrow\infty}\sin n\theta = -\frac{\theta^2}{2!}+\frac{\theta^4}{4!}- \frac{\theta^6}{6!}+...=\cos\theta - 1 = -2\sin\frac{1}{2}\theta
$$
and
\begin{equation}
\lim_{n\rightarrow\infty}\sin n\theta = -\frac{2\sin^{2}\frac{1}{2}\theta}{\theta}, \;\;\;\;\;\;\;\;\; (-\pi<\theta<\pi). \label{eq:trig_6}
\end{equation}

Analogously, taking the expansion of $\cos x$ into power series
$$
\cos x\theta = 1-\frac{x^2 \theta^2}{2!}+\frac{x^4 \theta^4}{4!}-...,
$$
where $0\leq \theta \leq 2\pi$, we get
\begin{equation}
\lim_{n\rightarrow\infty}\cos n\theta = \frac{2\sin\frac{1}{2}\theta \cos\frac{1}{2}\theta}{\theta}, \;\;\;\;\;\;\;\;\; (-\pi<\theta<\pi). \label{eq:trig_7}
\end{equation}

According to item (\ref{item:2}), for any real number $\theta$
$$
\sin(n+\beta)\theta = \sin n\theta \cos \beta\theta + \cos n\theta \sin \beta\theta.
$$
Then, passing to the limit and using (\ref{eq:trig_6}) and (\ref{eq:trig_7}), we have
\begin{eqnarray}
\lim_{n\rightarrow\infty}\sin(n+\beta)\theta = \cos \beta\theta \biggl(\frac{-2\sin^{2}\frac{1}{2}\theta}{\theta}\biggr) + \sin \beta\theta \biggl(\frac{2\sin\frac{1}{2}\theta \cos\frac{1}{2}\theta}{\theta}\biggr) =
 \frac{2\sin\frac{1}{2}\theta \sin(\beta-\frac{1}{2})\theta}{\theta} \nonumber
\end{eqnarray}
and
\begin{equation}
\lim_{n\rightarrow\infty}\sin(n+\beta)\theta = \frac{2\sin\frac{1}{2}\theta \sin(\beta-\frac{1}{2})\theta}{\theta}, \;\;\;\;\;\;\; (-\pi<\theta<\pi). \label{eq:trig_8}
\end{equation}

Analogously, relying on item (\ref{item:5}), we obtain the equality
\begin{equation}
\lim_{n\rightarrow\infty}\cos(n+\beta)\theta = \frac{2\sin\frac{1}{2}\theta \cos(\beta-\frac{1}{2})\theta}{\theta}, \;\;\;\;\;\;\; (-\pi<\theta<\pi). \label{eq:trig_9}
\end{equation}

From (\ref{eq:trig_8}) and (\ref{eq:trig_9}) the following equalities immediately follow
$$
\lim_{n\rightarrow\infty}\sin(n+\beta+1)\theta = -\lim_{n\rightarrow\infty}\sin(n-\beta)\theta, 
$$
$$
\lim_{n\rightarrow\infty}\cos(n+\beta+1)\theta = \lim_{n\rightarrow\infty}\cos(n-\beta)\theta
$$
and
$$
\lim_{n\rightarrow\infty}\sin(n+\frac{1}{2})\theta = 0,
$$
which are in good agreement with Lemma \ref{lemma2} and formula (\ref{eq:3_cor_lemma2}).

Substituting the formulas (\ref{eq:trig_6}), (\ref{eq:trig_7}), and (\ref{eq:trig_8}) into (\ref{eq:trig_1}) and (\ref{eq:trig_2}), we get the following two identities
$$
\frac{-2\sin\frac{1}{2}\theta \sin\theta}{\theta} \equiv \frac{-\sin\frac{1}{2}\theta (2\sin\frac{1}{2}\theta \cos\frac{1}{2}\theta)}{\theta},
$$
i. e.
$$
\frac{-2\sin\frac{1}{2}\theta \sin\theta}{\theta} \equiv \frac{-2\sin\frac{1}{2}\theta \sin\theta}{\theta}
$$
and
$$
\frac{-2\sin^{2}\frac{1}{2}\theta}{\theta} = -\tan\frac{1}{2}\theta \: \frac{2\sin\frac{1}{2}\theta \cos\frac{1}{2}\theta}{\theta},
$$
i. e.
$$
\frac{-2\sin^{2}\frac{1}{2}\theta}{\theta} \equiv \frac{-2\sin^{2}\frac{1}{2}\theta}{\theta},
$$
which confirm the equalities (\ref{eq:trig_1}) and (\ref{eq:trig_2}).

Using the formulas (\ref{eq:trig_6}) and (\ref{eq:trig_7}, that enable us to find the limits of $\sin n\theta$ and $\cos n\theta$, 
we obtain for the formulas (\ref{eq:trig_3}) and (\ref{eq:trig_4})
$$
\lim_{n\rightarrow\infty}\sin^{2}n\theta = \frac{1}{2}\cot\theta \Bigl(-\frac{\sin^{2}\theta}{\theta}\Bigr)+\frac{1}{2},
$$
i. e.
\begin{equation}
\lim_{n\rightarrow\infty}\sin^{2}n\theta = \frac{1}{2} - \frac{\sin 2\theta}{4\theta}, \;\;\;\;\;\;\;\;\; (-\pi<\theta<\pi) \label{eq:trig_10}
\end{equation}
and
\begin{equation}
\sum_{u=1}^{\infty}\sin u\theta = \frac{1}{2}\cot\frac{1}{2}\theta - \frac{1}{\theta}, \;\;\;\;\;\;\;\;\; (0<\theta<\pi). \label{eq:trig_11}
\end{equation}

Then, substituting (\ref{eq:trig_10}) into the formula (\ref{eq:trig_5}), we get
$$
\sum_{u=1}^{\infty}\sin(2u-1)\theta = \frac{1}{\sin\theta}\biggl(\frac{1}{2} - \frac{\sin\theta \cos\theta}{2\theta}\biggr),
$$
i. e.
$$
\sum_{u=1}^{\infty}\sin(2u-1)\theta = \frac{1}{2\sin\theta} - \frac{\cos\theta}{2\theta} \;\;\;\;\;\;\;\;\; (0<\theta<\pi).
$$

\begin{remark}
The formulas (\ref{eq:trig_6}), (\ref{eq:trig_7}), (\ref{eq:trig_10}) and others can be also obtained by direct integration, like in case of polynomials. Indeed,
$$
\int\limits_{-1}^{0}\sin\theta xdx = -\frac{1}{\theta}\cos\theta x \: \Big|_{-1}^{0} = \frac{-1}{\theta} (1-\cos\theta) = \frac{-2\sin^{2}\frac{1}{2}\theta}{\theta},
$$
$$
\int\limits_{-1}^{0}\cos\theta xdx = \frac{1}{\theta}\sin\theta x \: \Big|_{-1}^{0} = \frac{\sin\theta}{\theta}
$$
where $-\pi<\theta<\pi$.

And also
$$
\int\limits_{-1}^{0}\sin^{2}\theta xdx = \Bigl(\frac{x}{2}-\frac{1}{4\theta}\sin 2\theta x \Bigr) \: \Big|_{-1}^{0} = \frac{1}{2}-\frac{\sin 2\theta}{4\theta},
$$
where $-\frac{\pi}{2}<\theta<\frac{\pi}{2}$.
\end{remark}
\vspace{1cm}
\begin{center}
\rule{9.5cm}{0.25mm}
\end{center}

\newpage

\section{Extension of the Class of Regular Functions}\label{sec:exten}

In this section we give a way of how the class of regular functions can be extended when constructing them from non-elementary functions. In case of elementary functions, in Section \ref{sec:intro}, we construct regular functions $f(x)$, using the condition
\begin{equation}
F(z+1)-F(z)=f(z) \;\;\;\;\;\;\;\;\;\;\;\ \forall z\in \mathbb{Z}    \label{eq:regucond}
\end{equation}
where $F(x)$ is elementary.

To extend the class of regular functions by constructing them with use of non-elementary functions, we re-define the notion of regularity and the notion of odd/even function in the following way.
\begin{definition} \label{def:exten}
The function $f(x)$, $x\in \mathbb{Z}$, is called \textit{regular} if there exists a function $F(x)$, not necessarily elementary, such that 
$$
F(z+1)-F(z)=f(z) \;\;\;\;\;\;\;\;\; \forall z\in \mathbb{Z}
$$
and
$$
\lim_{n\rightarrow\infty}F(n+1)-\lim_{n\rightarrow\infty}F(n)=\lim_{n\rightarrow\infty}f(n).
$$
\end{definition}

\begin{definition}
A function $f(x)$, $x\in \mathbb{Z}$, is called \textit{even} if 
$$
f(-z)=f(z) \;\;\;\;\;\;\;\;\; \forall z\in \mathbb{Z}
$$
and
$$
\lim_{n\rightarrow\infty}f(-n)=\lim_{n\rightarrow\infty}f(n).
$$
\end{definition}

\begin{definition}
A function $f(x)$, $x\in \mathbb{Z}$, is called \textit{odd} if 
$$
f(-z)=-f(z) \;\;\;\;\;\;\;\;\; \forall z\in \mathbb{Z}
$$
and
$$
\lim_{n\rightarrow\infty}f(-n)=-\lim_{n\rightarrow\infty}f(n).
$$
\end{definition}

It is essential that in the case of elementary generating function, for functions $f(x)$ and $F(x)$ connected by (\ref{eq:regucond}) the limiting relation 
$$
\lim_{n\rightarrow\infty}F(n+1)-\lim_{n\rightarrow\infty}F(n)=\lim_{n\rightarrow\infty}f(n).
$$
is automatically satisfied. For this reason, it was not included in the definition of regular function given in Section \ref{sec:intro}.

So, we have the following
\begin{proposition} \label{prop:exten}
For any two elementary functions $f(x)$ and $F(x)$ such that \;
$F(z+1)-F(z)=f(z) \;\;\; \forall z\in \mathbb{Z}$
the following equality holds
\begin{equation}
\lim_{n\rightarrow\infty}F(n+1)-\lim_{n\rightarrow\infty}F(n)=\lim_{n\rightarrow\infty}f(n). \label{eq:prop_exten_main}
\end{equation}
\end{proposition}
\begin{proof}
This proposition will be proved by two ways.
\begin{itemize}
\item[(1)] In view of axiom (\ref{axiom:5}), we have
\begin{equation}
F(n+1)-F(1)=\sum_{u=1}^{n}f(u) \;\;\;\;\;\;\;\;\; \forall n  \label{eq:prop_exten_1}
\end{equation}
and
\begin{equation}
F(n)-F(1)=\sum_{u=1}^{n}f(u)-f(n) \;\;\;\;\;\;\;\;\; \forall n. \label{eq:prop_exten_2}
\end{equation}
Subtracting (\ref{eq:prop_exten_2}) from (\ref{eq:prop_exten_1}), we get
$$
S_n = F(n+1)-F(n)-f(n)=\sum_{u=1}^{n}\omega(u) \;\;\;\;\;\;\; \forall n,
$$
where $\omega(u)\equiv 0$.

Now passing to the limit, we have
$$
\lim_{n\rightarrow\infty}S_n=\lim_{n\rightarrow\infty}F(n+1)-\lim_{n\rightarrow\infty}F(n) - \lim_{n\rightarrow\infty}f(n)=\sum_{u=1}^{\infty}\omega(u).
$$
And since $\sum\limits_{u=1}^{\infty}\omega(u)=0$, we finally obtain
$$
\lim_{n\rightarrow\infty}F(n+1)-\lim_{n\rightarrow\infty}F(n)=\lim_{n\rightarrow\infty}f(n),
$$
which completes the proof. \\
\item[(2)] According to axiom (\ref{axiom:5}) and formula (\ref{eq:1_from_lemma3}), we have
\begin{equation}
F(-n)-F(1)=\sum_{u=1}^{-(n+1)}f(u)=-\sum_{u=0}^{n}f(-u) \;\;\;\;\;\;\;\;\; \forall n  \label{eq:prop_exten_3}
\end{equation}
and
\begin{equation}
F(n)-F(1)=\sum_{u=1}^{n}f(u)-f(n) \;\;\;\;\;\;\;\;\; \forall n. \label{eq:prop_exten_4}
\end{equation}
Subtracting (\ref{eq:prop_exten_4}) from (\ref{eq:prop_exten_3}), we get
$$
F(-n)-F(n)-f(n)=-\biggl(\sum_{u=0}^{n}f(-u)+\sum_{u=1}^{n}f(u)\biggr) \;\;\;\;\;\;\;\;\; \forall n.
$$
Passing to the limit and using Lemma \ref{lemma2}, we obtain
$$
\lim_{n\rightarrow\infty}F(n+1)-\lim_{n\rightarrow\infty}F(n) - \lim_{n\rightarrow\infty}f(n)=
-\biggl(\sum_{u=1}^{\infty}f(-u)+\sum_{u=0}^{\infty}f(u)\biggr),
$$
which in view of (\ref{eq:1_from_lemma1}) finally gives
$$
\lim_{n\rightarrow\infty}F(n+1)-\lim_{n\rightarrow\infty}F(n)=\lim_{n\rightarrow\infty}f(n).
$$
This completes the proof.
\end{itemize}
The proposition is proved.
\end{proof}

However, in general case, for the case of non-elementary generating function, the condition (\ref{eq:prop_exten_main}) should be explicitly included in the definition of regular function as it is done in Definition \ref{def:exten}.

Thus, the class of regular functions can be substantially extended by using non-elementary functions.
\vspace{1cm}
\begin{center}
\rule{9.5cm}{0.25mm}
\end{center}

\newpage

\section{Conclusion}\label{sec:conc}

In this work we presented the fundamentals of the new field of research leading to a new theory that arises from the creating a new method for ordering the
integers \citep{E,P,bag,bagR,bagturkey}. This method, which relies on a single basic principle, allows one to obtain the limits of
unbounded and oscillating functions and to sum divergent series.

The tools and techniques
developed in the framework of our conception allowed us
to assert that every elementary function of integer argument given on $\mathbb{Z}$ has
a definite limit, to compute the limits of unbounded and oscillating functions,
to discover a unified and general method for summation of divergent series,
to find several properties of divergent series, and some other results.
The investigations within the framework
of this field reveal a number of unexpected 
phenomena concerning by now several areas of mathematics. 

It should be mentioned that the selection of topics in this work is of course
quite arbitrary and has the aim of demonstration of the efficiency of the new
method of ordering the integers, how it works, and its potential for
mathematical research. On this topics we also tried to show the consistency of the theory presented, 
in spite of seeming anomalousness of some results.

All the results of this work, even those for obtaining
of which as a rule complex-analytical techniques are employed, have been
obtained in elementary way, without addressing the notion of analytic continuation.
For example, in \citep{bag} we concerned with the values of the Riemann zeta function at negative integers,
but in our more comprehensive work \citep{bagR} we will deal with the values of the zeta function not only at negatives 
but at positive integers as well.

The work has an introductory character and the theoretical basics 
presented here is a good setting for future research. 

We have left open the question of the physical meaning of connection of negative and positive numbers 
through infinity; in fact, it may be the question of proof of space closedness, or perhaps 
it concerns the dual essence of the physical space. The study and analysis of physical aspects
we leave for future works.
However, if one would assume that our axiom system is in accordance with
the real physical world, the equalities obtained in Lemmas \ref{lemma:1} and \ref{lemma2}
could serve as an evidence that the physical space, being unbounded, is closed and finite.
It is worth mentioning here that the analogous idea was expressed by Riemann in 1854 in his lecture
``On the hypotheses underlying geometry'' in which he supposed that the space, albeit unbounded, is finite, and the line
acquires a finite length. The importance of these ideas was fully revealed later, 
after Einstein made use of them to establish the foundations of his general relativity theory.

The techniques we have developed, and further refinements of them, can also be used to study double, miltiple, 
and functional series. Another direction, on which we paid attention in the work, is the theory of infinite
products. The ideas and results of the work give hope to shed light
on and obtain a deeper understanding of the structure of numerical continuum. 
From the obtained formulas for Bernoulli numbers other identities and congruences with Bernoulli 
and Euler numbers can be derived \citep{Lan}.

It would be of interest to study the topological properties of  new
number axis, to construct geometric models based on it, and to study 
geometric objects in the light of this theory. Then, it would be interesting
to study how this theory and other branches of mathematics, especially,
complex analysis relate to each other. This direction is interesting in itself,
especially because  our results confirm some results obtained in the theory of the Riemann's sphere \citep{J}.

We hope that this work will be of certain interest for mathematical community, and 
further development of the theory will have a wide spectrum of directions.


We believe that further development of the theory will not just provide
new interesting problems within theory itself but will also have an impact
on other domains of mathematics.

\ack The second author is very indebted to R. Varshamov for inviting and involving him in this new field of research.


\end{document}